\documentclass[letterpaper,11pt,oneside,
              ]{amsart}

\usepackage{amssymb}
\usepackage{amsxtra}
\usepackage[mathscr]{eucal}
\usepackage{graphics}

\makeatletter
\renewcommand\@setcopyright{}
\renewcommand{\tocsection}[3]{%
  \indentlabel{\@ifnotempty{#2}{\ignorespaces#1 #2.\quad}}#3\dotfill}
\makeatother

\numberwithin{equation}{section}

\newtheorem{theorem}{Theorem}[section]
\newtheorem{proposition}[theorem]{Proposition}
\newtheorem{lemma}[theorem]{Lemma}

\theoremstyle{definition}

\newtheorem{definition}[theorem]{Definition}
\newtheorem{example}[theorem]{Example}

\newtheorem{remark}[theorem]{Remark}
\newtheorem{remarks}[theorem]{Remarks}

\newcommand{\al}{\alpha}
\newcommand{\aM}{\alpha_M}
\newcommand{\ann}{\operatorname{ann}}
\newcommand{\arda}{\al_{\rd/\fa}}
\newcommand{\ardfm}{\al_{\rdfm}}
\newcommand{\ardp}{\al_{\rd/\fp}}
\newcommand{\arof}{\al_{R_1/\langle f\rangle}}
\newcommand{\artf}{\al_{R_2/\langle f\rangle}}
\newcommand{\artp}{\al_{R_2/\fp}}
\newcommand{\arthp}{\al_{R_3/\fp}}
\newcommand{\arthf}{\al_{R_3/\langle f\rangle}}
\newcommand{\ardf}{\al_{R_d/\langle f\rangle}}
\newcommand{\asc}{\operatorname{asc}}
\newcommand{\be}{\mathbf{e}}
\renewcommand{\bf}{\mathbf{f}}
\newcommand{\bk}{\mathbf{k}}
\newcommand{\bm}{\mathbf{m}}
\newcommand{\bn}{\mathbf{n}}
\newcommand{\bp}{\mathbf{p}}
\newcommand{\bs}{\mathbf{s}}

\newcommand{\bx}{\mathbf{x}}
\newcommand{\bv}{\mathbf{v}}
\newcommand{\bw}{\mathbf{w}}
\newcommand{\bz}{\mathbf{z}}
\newcommand{\CC}{\mathbb{C}}
\newcommand{\charact}{\operatorname{char}}
\newcommand{\D}{\mathsf{D}}
\newcommand{\Da}{\Delta_{\alpha}}
\newcommand{\Db}{\Delta_{\beta}}
\newcommand{\diam}{\operatorname{diam}}
\newcommand{\dist}{\operatorname{dist}}
\newcommand{\entropy}{\textsl{\sffamily h}}
\newcommand{\entrk}{\operatorname{entrk}}

\newcommand{\exprk}{\operatorname{exprk}}
\newcommand{\E}{\mathsf{E}}
\newcommand{\fa}{\mathfrak{a}}
\newcommand{\fb}{\mathfrak{b}}
\newcommand{\ff}{\mathfrak{f}}
\newcommand{\fp}{\mathfrak{p}}
\newcommand{\fq}{\mathfrak{q}}
\newcommand{\ftil}{\widetilde{f}}
\newcommand{\G}{\mathsf{G}}
\renewcommand{\H}{\mathsf{H}}
\newcommand{\Hd}{\mathsf{H}_d}
\newcommand{\Hv}{H_{\bv}}
\newcommand{\Hw}{H_{\bw}}
\newcommand{\hz}{H_{\mathbb{Z}}}
\def\implies{ \ensuremath{\Rightarrow} }
\newcommand{\ifff}{\ensuremath{\Leftrightarrow}}
\newcommand{\kdim}{\operatorname{kdim}}
\newcommand{\kernel}{\operatorname{ker}}

\newcommand{\lt}{\operatorname{lt}}
\newcommand{\N}{\mathsf{N}}
\newcommand{\Nn}{\mathsf{N}^{\mathsf{n}}}
\newcommand{\Nv}{\mathsf{N}^{\mathsf{v}}}
\newcommand{\NN}{\mathbb{N}}
\newcommand{\Newt}{\mathscr{N}}
\newcommand{\pg}{{\fp}_{\Gamma}}
\newcommand{\rd}{R_d}
\newcommand{\rda}{R_d/\fa}
\newcommand{\rdf}{{R_d/\langle f\rangle}}
\newcommand{\rdfm}{{R_d/\ff(M)}}
\newcommand{\rdp}{R_d/\fp}
\newcommand{\rthp}{R_3/\fp}
\newcommand{\QQ}{\mathbb{Q}}
\renewcommand{\Re}{\operatorname{Re}}
\newcommand{\rg}{R_{\Gamma}}
\newcommand{\rh}{R_H}
\newcommand{\RR}{\mathbb R}
\newcommand{\sph}{\mathsf{S}}
\newcommand{\sphd}{\mathsf{S}_{d-1}}
\newcommand{\SA}{\mathscr{A}}
\newcommand{\SC}{\mathscr{C}}
\newcommand{\SCQ}{\SC_{\QQ}}

\newcommand{\SN}{\mathscr{N}}
\def\SS{\mathbb{S}}   
\newcommand{\SU}{\mathscr{U}}
\newcommand{\supp}{\operatorname{supp}}
\newcommand{\SV}{\mathscr{V}}
\newcommand{\SW}{\mathscr{W}}
\newcommand{\TT}{\mathbb{T}}

\newcommand{\V}{\mathsf{V}}
\newcommand{\veps}{\varepsilon}
\newcommand{\vH}{\mathbf{v}_{\!H}^{}}

\newcommand{\xbar}{\overline{x}}
\newcommand{\xtri}{x^{\vartriangle}}
\newcommand{\ybar}{\overline{y}}
\newcommand{\zd}{\mathbb{Z}^d}
\newcommand{\ZZ}{\mathbb{Z}}
\newcommand{\ztt}{\ZZ[t,1/2t(t+1)]}

\newcommand{\<}{\langle}
\renewcommand{\>}{\rangle}  

\newsymbol\emptyphi 203F
\renewcommand{\emptyset}{\emptyphi}

\begin{document}

\allowdisplaybreaks

\title[Expansive subdynamics for algebraic actions]
  {Expansive subdynamics for algebraic
  $\mathbb Z^d$-actions}

\author{Manfred Einsiedler, Douglas Lind, Richard Miles, and
  Thomas Ward}

\address{Manfred Einsiedler, Mathematisches Institut,
  Universit\"at Wien, Strudlhofgasse 4, A-1090, Vienna, Austria}
\email{manfred@mat.univie.ac.at}

\address{Douglas Lind, Department of Mathematics,
Box 354350, University of Washington, Seattle, WA 98195--4350, USA}
\email{lind@math.washington.edu}

\address{Richard Miles, School of Mathematics, University of East
  Anglia, Norwich NR4 7TJ, United Kingdom}
\email{r.miles@uea.ac.uk}

\address{Thomas Ward, School of Mathematics, University of East
  Anglia, Norwich NR4 7TJ, United Kingdom}
\email{t.ward@uea.ac.uk}

\thanks{The first author was supported in part by FWF Research
Grant P12250--MAT and EPSRC Grant GR/M 49588. The second author
was supported in part by NSF Grant DMS-9622866.  The first two
authors gratefully acknowledge support from the Schr\"odinger
Institute at the University of Vienna.  The third author
gratefully acknowledges support from EPSRC Grant 97700813. The
fourth author thanks the University of East Anglia and the
American Mathematical Society for their support.}

\date{\today}



\begin{abstract}
   A general framework for investigating topological actions of
   $\ZZ^d$ on compact metric spaces was proposed by Boyle and
   Lind in terms of expansive behavior along lower-dimensional
   subspaces of $\RR^d$. Here we completely describe this
   expansive behavior for the class of algebraic $\ZZ^d$-actions
   given by commuting automorphisms of compact abelian groups.
   The description uses the logarithmic image of an algebraic
   variety together with a directional version of Noetherian
   modules over the ring of Laurent polynomials in several
   commuting variables.
   
   We introduce two notions of rank for topological
   $\zd$-actions, and for algebraic $\zd$-actions describe how
   they are related to each other and to Krull dimension. For a
   linear subspace of $\RR^d$ we define the group of points
   homoclinic to zero along the subspace, and prove that this
   group is constant within an expansive component.
\end{abstract}

\maketitle

{\setlength{\linewidth}{2in}
\tableofcontents}

\section{Introduction}\label{sec:introduction}

Expansiveness is a multifaceted property that plays an
important role throughout dynamics. Let $\beta$ be an action of
$\zd$ by homeomorphisms of a compact metric space $(X,\rho)$.
Then $\beta$ is called \textit{expansive} if there is a
$\delta>0$ such that if $x$ and $y$ are distinct points of $X$
then there is an $\bn\in\zd$ such that $\rho(\beta^\bn x,\beta^\bn
y)>\delta$.

In \cite{BL} the notion of expansiveness along a subset, and
especially a subspace, of $\RR^d$ was introduced, by considering
only those elements of $\zd$ that lie within a given bounded distance
of the set. Let $\G_k$ denote the compact Grassmann manifold of
$k$-dimensional subspaces (or \textit{$k$-planes}) in $\RR^d$,
and $\N_k(\beta)$ be the set of those $k$-planes which are
not expansive for $\beta$. It was shown in \cite{BL} that if $X$
is infinite, then $\N_{d-1}(\beta)$ is a nonempty compact subset
of $\G_{d-1}$ that determines all other $\N_k(\beta)$ as follows:
a $k$-plane is nonexpansive for $\beta$ if and only if it is
contained in some subspace in ~$\N_{d-1}(\beta)$.

Denote by $\E_k(\beta)$ the set of expansive $k$-planes for
$\beta$, which is an open subset of~ $\G_k$. Various dynamical
notions, such as entropy, can be defined along subspaces. The
expansive subdynamics philosophy advocated in \cite{BL} proposes
that many such properties should be either constant or vary
nicely within a connected component of ~$\E_k(\beta)$, but that
they should typically change abruptly when passing from one
component to another, analogous to a ``phase transition.''
Several examples of this philosophy in action are given in
\cite{BL}.  In Section~ \ref{sec:homoclinic-group} we provide
another by considering points homoclinic along subspaces.  Thus a
basic starting point in the analysis of any topological
$\ZZ^d$-action is to describe its expansive subspaces, especially
those of co-dimension one, since these determine the rest.

An \textit{algebraic $\zd$-action} is an action $\alpha$ of $\zd$
by (continuous) automorphisms of a compact abelian group, which
we assume to be metrizable. We will consistently use $\alpha$ to
denote an algebraic $\zd$-action and $\beta$ for a general
topological $\zd$-action. Such algebraic actions have provided a
rich source of examples and phenomena (see \cite{Sch}). The
purpose of this paper is to completely determine the expansive
subspaces for all algebraic $\zd$-actions.

In Section \ref{sec:expansive-subdynamics} we review the relevant
ideas from \cite{BL}, and show that it is sufficient to use
half-spaces rather than $(d-1)$-dimensional planes. The algebra
we need is described in Section~ \ref{sec:algebraic-actions}. In
Section~ \ref{sec:characterization} we develop our main result,
Theorem~\ref{thm:characterization}, which describes expansive
half-spaces in terms of prime ideals.  We give in Section~
\ref{sec:examples} a number of examples that illustrate and
motivate our results from the previous section. In Section~
\ref{sec:further-analysis} we investigate the prime ideal case in
more detail, including an algorithm to compute the expansive set.
We introduce two notions of ``rank'' for a topological
$\zd$-action in Section \ref{sec:ranks}, and for algebraic
$\zd$-actions show how they are related to each other and to
Krull dimension. In Section \ref{sec:lower-dimensional} we
extend our basic results to lower-dimensional subspaces of
$\RR^d$.  The homoclinic group along a subspace is defined in
Section~ \ref{sec:homoclinic-group} and shown to be constant
within an expansive component. This fact has some interesting
dynamical consequences.

\section{Expansive subdynamics}\label{sec:expansive-subdynamics}

Let $(X,\rho)$ be a compact metric space, which we assume is
infinite unless otherwise stated. A \textit{$\zd$-action} $\beta$
on $X$ is a homomorphism from $\zd$ to the group of
homeomorphisms of $X$. For $\bn\in\zd$ we denote the
corresponding homeomorphism by $\beta^{\bn}$, so that
$\beta^{\bm}\circ\beta^{\bn}=\beta^{\bm+\bn}$, and
$\beta^{\mathbf{0}}$ is the identity on~ $X$. For a subset $F$ of
$\RR^d$ put
\begin{displaymath}
   \rho_\beta^F(x,y)=\sup\{\rho\bigl(\beta^{\bn}(x),\beta^{\bn}(y)
   \bigr):
   \bn\in F\cap\zd\},
\end{displaymath}
and if $F\cap\zd=\emptyset$ define $\rho_\beta^F(x,y)=0$.

\begin{definition}\label{def:expansive}
   A $\zd$-action $\beta$ on $(X,\rho)$ is \textit{expansive}
   provided there is a $\delta>0$ such that
   $\rho_\beta^{\RR^d}(x,y)\le\delta$ implies that $x=y$. In this
   case $\delta$ is called an \textit{expansive constant} for ~$\beta$.
\end{definition}

Let $\|\cdot\|$ denote the Euclidean norm on $\RR^d$. For
$F\subset \RR^d$ and $\bx\in\RR^d$ define
$\dist(\bx,F)=\inf\{\|\bx-\mathbf{y}\|:\mathbf{y}\in F\}$.
For $t>0$ put $F^t=\{\bx\in\RR^d:\text{dist}(\bx,F)\le t\}$, so
that $F^t$ is the result of thickening $F$ by ~$t$.

\begin{definition}\label{def:expansive-subset}
   Let $\beta$ be a $\zd$-action on $(X,\rho)$ and $F$ be a
   subset of $\RR^d$. Then $F$ is \textit{expansive for} $\beta$,
   or $\beta$ is \textit{expansive along} $F$, if there are
   $\varepsilon>0$ and $t>0$ such that
   $\rho_\beta^{F^t}(x,y)\le\varepsilon$ implies that $x=y$. If $F$
   fails to meet this condition it is 
   \textit{nonexpansive for} $\beta$, or $\beta$ is
   \textit{nonexpansive along ~$F$}.
\end{definition}

\begin{remark}\label{rem:subset}
   Every subset of a nonexpansive set for $\beta$ is clearly also
   nonexpansive for ~$\beta$.  Every translate of an expansive set
   is expansive \cite[p.\ 57]{BL}. In the above definition we can
   take for $\varepsilon$ a fixed expansive constant for ~$\beta$
   \cite[Lemma~ 2.3]{BL}.
\end{remark}

Next we examine subsets $F$ that are linear subspaces of $\RR^d$.
Let $\G_k=\G_{d,k}$ denote the Grassmann manifold of
$k$-dimensional subspaces (or simply \textit{$k$-planes}) of
$\RR^d$. Recall that $\G_k$ is a compact manifold of dimension
$k(d-k)$ whose topology is given by declaring two subspaces to be
close if their intersections with the unit sphere are close in the
Hausdorff metric. A $k$-plane and its $(d-k)$-dimensional
orthogonal complement determine each other, giving a natural
homeomorphism between $\G_k$ and ~$\G_{d-k}$.

\begin{definition}\label{def:expansive-set}
   For a $\zd$-action $\beta$ define
   \begin{align*}
      \E_k(\beta)&=\{V\in\G_k:\text{$V$ is expansive
        for $\beta$}\},\\
      \N_k(\beta)&=\{V\in\G_k:\text{$V$ is nonexpansive
        for $\beta$}\}.
   \end{align*}
   An \textit{expansive component of $k$-planes for $\beta$} is a
   connected component of $\E_k(\beta)$.
\end{definition}

\begin{example}\label{exam:ledrappier}
   (\textit{Ledrappier's example}.)
   Take $d=2$,
   \begin{displaymath}
      X=\{x\in(\ZZ/2\ZZ)^{\ZZ^2}:x_{i,j}+x_{i+1,j}+x_{i,j+1}\equiv
      0\text{\, (mod 2) for all $i,j$}\},
   \end{displaymath}
   and let $\beta$ be the $\ZZ^2$-action generated by the horizontal
   and vertical shifts. If $L$ is a line that is not parallel to
   one of the sides of the unit simplex in $\RR^2$ and $t\ge2$,
   then for each $x\in X$ the coordinates of $x$ within $L^t$
   determine all of ~$x$, so that $L\in\E_1(\beta)$. On the other
   hand, the three lines parallel to the sides of the simplex do
   not have this property, and they comprise $\N_1(\beta)$ (see
   \cite[Example~ 2.7]{BL} for details).
\end{example}

Simple coding arguments \cite[Lemma~ 3.4]{BL} show that each
$\E_k(\beta)$ is an open subset of ~$\G_k$, so that each
$\N_k(\beta)$ is compact.  Hence expansive components of
$k$-planes for $\beta$ are open subsets of $\G_k$.  By Remark
\ref{rem:subset}, if $W$ is nonexpansive for ~$\beta$ and $V$ is
a subspace of ~$W$, then $V$ is also nonexpansive for ~$\beta$. A
basic result \cite[Theorem\ 3.6]{BL} is the converse: If $V$ is a
nonexpansive subspace for~ $\beta$ of dimension $\le d-2$, then
there is a nonexpansive subspace for $\beta$ containing $V$ of
one higher dimension. If $X$ is infinite, then the zero subspace
is nonexpansive, and hence inductively we see that each
$\N_k(\beta)\ne\emptyset$ for $1\le k\le d-1$. Furthermore, it
follows that if $V\in\N_k(\beta)$, then there is a
$W\in\N_{d-1}(\beta)$ that contains ~$V$. Hence $\N_k(\beta)$
consists of exactly all $k$-dimensional subspaces of the
subspaces in ~$\N_{d-1}(\beta)$.  Thus $\N_{d-1}(\beta)$
determines the entire expansive subdynamics of $\beta$.

In order to treat algebraic $\zd$-actions, it is convenient to
shift our viewpoint slightly and use half-spaces in $\RR^d$
rather than $(d-1)$-planes. Let $\sphd=\{\bv\in\RR^d:\|\bv\|=1\}$
be the unit $(d-1)$-sphere. For $\bv\in\sphd$ define
$\Hv=\{\bx\in\RR^d: \bx\cdot\bv\le0\}$ to be the half-space with
outward unit normal ~$\bv$. Let $\Hd$ be the set of half-spaces in
~$\RR^d$, which we identify with $\sphd$ via the parameterization
$\bv \leftrightarrow\Hv$. For $H\in\Hd$ we denote its outward
unit normal vector by ~$\vH$.

Expansiveness along a half-space $H$ is defined using Definition
\ref{def:expansive-set} with $F=H$. Observe that thickening $\Hv$
by $t>0$ results merely in the translation $\Hv+t\bv$ of
$\Hv$. Hence there is no need to thicken half-spaces in the
definition, and a $\zd$-action $\beta$ is therefore expansive
along $H$ if and only if there is an $\varepsilon>0$ such that
$\rho_{\beta}^{H}(x,y)\le\varepsilon$ implies that $x=y$.

\begin{definition}\label{def:expansive-half-space}
   For a $\zd$-action $\beta$ define
   \begin{align*}
      \E(\beta)&=\{H\in\Hd:\text{$H$ is expansive
        for $\beta$}\},\\
      \N(\beta)&=\{H\in\Hd:\text{$H$ is nonexpansive
        for $\beta$}\}.
   \end{align*}
   An \textit{expansive component of half-spaces for $\beta$} is
   a connected component of $\E(\beta)$.
\end{definition}

\begin{remark}
  A coding argument analogous to \cite[Lemma~ 3.4]{BL} shows that
  $\E(\beta)$ is an open set and so $\N(\beta)$ is a compact set.
\end{remark}

The following lemma shows that a $(d-1)$-plane is nonexpansive
for ~$\beta$ if and only if at least one of the two bounding
half-spaces is also nonexpansive for ~$\beta$. Thus if we define
$\pi\colon \Hd\to\G_{d-1}$ by $\pi(H)=\partial H$, then
$\pi(\N(\beta))=\N_{d-1}(\beta)$. This shows that the half-space
behavior $\N(\beta)$ determines the expansive subdynamics of~
~$\beta$.

We start by recalling the following key notion from
\cite[Definition 3.1]{BL}.

\begin{definition}\label{def:coding}
   Let $\beta$ be an expansive $\zd$-action with expansive
   constant ~$\delta$. For subsets $E$, $F$ of $\RR^d$ we say that
   $E$ \textit{codes} $F$ provided that, for every $\bx\in\RR^d$, if
   $\rho_\beta^{E+\bx}(x,y)\le\delta$ then
   $\rho_\beta^{F+\bx}(x,y)\le\delta$.
\end{definition}

\begin{lemma}\label{lem:boundary-half-space}
   Let $\beta$ be a $\zd$-action and $V\in \G_{d-1}$. Then
   $V\in\N_{d-1}(\beta)$ if and only if there is an
   $H\in\N(\beta)$ with $\partial H=V$.
\end{lemma}

\begin{proof}
   If $H\in\N(\beta)$, then by Remark \ref{rem:subset} we see
   that $V=\partial H\subset H$ is also nonexpansive.

   Conversely, let $V\in\G_{d-1}$ and $H=\Hv$, $H'=H_{-\bv}$ be
   the two half-spaces with boundary ~$V$. Suppose that both $H$
   and $H'$ are expansive for ~$\beta$. We prove that $V$ is also
   expansive for ~$\beta$, which will complete the proof.
   
   Since $\beta$ has an expansive half-space, it is an expansive
   action.  Let $\delta>0$ be an expansive constant for ~$\beta$.
   Let $B(r)$ denote the ball of radius~ $r$ in ~$\RR^d$, and
   $[\mathbf{0},\bv]$ be the line segment joining $\mathbf{0}$ to
   ~$\bv$. A ``finite'' version of the expansiveness of $H$,
   entirely analogous to \cite[Lemma~ 3.2]{BL}, is that there is
   an $r>0$ such that $H\cap B(r)$ codes $[\mathbf{0},\bv]$.
   Similarly, there is an $s>0$ such that $H'\cap B(s)$ codes
   $[\mathbf{0},-\bv]$. Hence if $t=\max\{r,s\}$, then $V^t$
   codes ~$V^{t+1}$, which by the same argument codes ~$V^{t+2}$,
   and so on.  Thus $V^t$ codes ~$\RR^d$, which means that $V$ is
   expansive.
\end{proof}

\section{Algebraic $\zd$-actions}\label{sec:algebraic-actions}

An \textit{algebraic $\zd$-action} is an action of $\zd$ by
(continuous) automorphisms of a compact abelian group. Such
actions provide a rich class of examples of $\zd$-actions having
striking connections with commutative algebra. The monograph by
Schmidt \cite{Sch} provides a detailed account of this theory.
For background and standard results from
commutative algebra used below the reader may consult \cite{Eis}.

Let $X$ be a compact abelian group with identity element ~$0_X$.
Suppose that $\al$ is an action of $\zd$ by automorphisms of $X$.
Let $M=\widehat{X}$, the Pontryagin dual group of ~$X$.  Define
$\rd=\ZZ[u_1^{\pm1},\dots,u_d^{\pm1}]$, the ring of Laurent
polynomials with integer coefficients in the $d$ commuting
variables $u_1$, $\dots$, $u_d$. We can make $M$ into an
$\rd$-module by defining $u_j\cdot m=\widehat{\al}^{\be_j}(m)$
for all $m\in M$, where $\be_j\in\zd$ is the $j$th unit vector
and $\widehat{\al}^{\be_j}$ is the automorphism of $M$ dual to
~$\al^{\be_j}$. An element $f\in\rd$ has the form
$f=f(u)=\sum_{\bn\in\zd}c_f(\bn)u^{\bn}$, where the $c_f(\bn)\in\ZZ$
and $c_f(\bn)=0$ for all but finitely many $\bn\in\zd$, and
$u^{\bn}=u_1^{n_1}\dots u_d^{n_d}$. Then $f\cdot
m=\sum_{\bn\in\zd} c_f(\bn)\widehat{\al}^{\bn}(m)$ for every
$m\in M$.

This process can be reversed. Suppose that $M$ is an
$\rd$-module. Let $X_M=\widehat{M}$ be its compact abelian dual
group. Each $u_j$ is a unit in ~$\rd$, so the map $\gamma_j$
defined by $\gamma_j(m)=u_j\cdot m$ is an automorphism of ~$M$.
Define an algebraic $\zd$-action $\al_M$ on $X_M$ by
$\al_M^{\be_j}=\widehat{\gamma}_j$. See \cite[Chapter\ II]{Sch}
for further explanation and many examples.

Hence using duality we see there is a one-to-one correspondence
between algebraic $\zd$-actions on the one hand and $\rd$-modules
on the other.  

A module over an arbitrary ring is \textit{Noetherian} if it
satisfies the ascending chain condition for submodules.  The ring
$\rd$ is Noetherian as a module over itself. A prime ideal
$\fp\subset\rd$ is \textit{associated to~ $M$} if there is an
$m\in M$ with $\fp=\{f\in\rd:f\cdot m=0\}$.  Let $\asc(M)$ denote
the set of prime ideals associated to an $\rd$-module $M$.  If
$M$ is Noetherian, then $\asc(M)$ is finite. One basic discovery
has been that the dynamical properties of $\al_M$ can largely be
determined from $\asc(M)$.

For example, let us describe when $\al_M$ is expansive. In order
to do so, we need some notation. Let $\CC$ denote the complex
numbers and $\CC^\times=\CC\smallsetminus\{0\}$. For an ideal
$\fa$ in $\rd$ put
\begin{displaymath}
   \V(\fa)=\bigl\{\,\bz=(z_1,\dots,z_d)\in(\CC^\times)^d:f(z_1,\dots,z_d)=0
   \text{ for all $f\in\fa$}\,\bigr\}
\end{displaymath}
(we omit $0$ from $\CC$ since we are using Laurent polynomials).
Let $\SS=\{z\in\CC:|z|=1\}$ and
$\SS^d=\{(z_1,\dots,z_d)\in\CC^d:|z_1|=\dots=|z_d|=1\}$ be the
multiplicative $d$-torus. Then we have the following
characterization of expansiveness due to Schmidt (see
\cite[Theorem~6.5]{Sch}), which can be thought as a
generalization of the fact that a toral automorphism is expansive
if and only if the associated integer matrix has no eigenvalue in~ $\SS$.

\begin{theorem}\label{thm:expansive}
   Let $M$ be an $\rd$-module. Then the following are equivalent.
   \begin{enumerate}
     \item[(1)] $\al_M$ is expansive.
     \item[(2)] $M$ is Noetherian and $\al_{\rd/\fp}$ is expansive for
      every $\fp\in\asc(M)$.
     \item[(3)] $M$ is Noetherian and $\V(\fp)\cap\SS^d=\emptyset$ for
      every $\fp\in\asc(M)$. 
   \end{enumerate}
\end{theorem}

Since the expansive subdynamics of a nonexpansive action are
trivial, we will consider only Noetherian $\rd$-modules from now
on. In particular, a Noetherian $\rd$-module $M$ is countable, and
hence $X_M$ is metrizable. We fix a metric $\rho$ on $X_M$,
which we may assume to be translation-invariant.
Observe that
$\rho\bigl(\al^{\bn}(x),\al^{\bn}(y)\bigr)=\rho\bigl(
\al^{\bn}(x-y),\al^{\bn}(0_X)\bigr)$ by our assumption on~ $\rho$.
Thus in considering the expansive
behavior of an algebraic $\zd$-action on a pair of points, we may
assume that one of them is ~$0_X$.

Roughly speaking, nonexpansiveness of $\alpha_M$ can occur for
two reasons: algebraic and geometric. The algebraic reason occurs
when $M$ is not Noetherian: in this case there is a decreasing sequence
of closed, $\alpha_M$-invariant subgroups converging to
$\{0_X\}$, and this immediately provides, for every $\delta>0$, a
nonzero point that remains within $\delta$ of $0_X$ under all
iterates of $\alpha_M$. The geometric reason occurs when
$\V(\fp)\cap\SS^d\ne\emptyset$: in this case it is possible to use a
element from $\V(\fp)\cap\SS^d$ to construct points that remain
arbitrarily close to $0_X$ under all iterates (see the proof of
Theorem \ref{thm:characterization} for details). For a
valuation-theoretic approach to expansive behavior see
\cite{Miles}.

\section{Characterization of expansive half-spaces}
   \label{sec:characterization}

Let $M$ be an $\rd$-module and $\al_M$ be the corresponding
algebraic $\zd$-action. In this section we characterize those
half-spaces $H\in\Hd$ which are expansive for $\al_M$ in terms of the
prime ideals associated to ~$M$. The main result, Theorem~
\ref{thm:characterization}, is a ``one-sided'' version of
Theorem~ \ref{thm:expansive}. The reader is urged to consult the
examples in Section~ \ref{sec:examples} first to motivate what
follows in this section.

According to Theorem~ \ref{thm:expansive}, there are two reasons
that $\al_M$ may fail to be expansive: $M$ may not be
Noetherian, or there may be a point in $\V(\fp)\cap\SS^d$ for
some $\fp\in\asc(M)$. In the following sequence of results we
investigate a ``one-sided'' version of each of these
possibilities. Our proofs closely parallel those in \cite{Sch},
with suitable modifications for their one-sided nature. 

We start with the Noetherian condition. For $H\in\Hd$, recall that
$\vH$ denotes the outward unit normal for ~$H$. Define
$\hz=H\cap\zd$.  Put $\rh=\ZZ[u^{\bn}:\bn\in\hz]$, which is a
subring of ~$\rd$.  It is important to note that in general $\rh$
is \textit{not} a Noetherian ring. Indeed, $\rh$ is Noetherian
exactly when $\vH$ is a rational direction in the sense that
$\RR\vH\cap\zd\ne\{\mathbf{0}\}$, so that $\rh$ is Noetherian
for only countably many ~$H$. Understanding when an
$\rd$-module is Noetherian as an $\rh$-module (i.e., when it is
\textit{$\rh$-Noetherian}) is one of the key points in our analysis.

\begin{lemma}\label{lem:Rh-noetherian}
   Let $M$ be a Noetherian $\rd$-module and $H\in\Hd$. Then $M$
   is $\rh$-Noetherian if and only if $\rd/\fp$ is
   $\rh$-Noetherian for every $\fp\in\asc(M)$.
\end{lemma}
 
\begin{proof}
   Suppose that $M$ is $\rh$-Noetherian. If $\fp\in\asc(M)$, then
   there is an $m\in M$ such that $\rd\cdot m\cong \rd/\fp$. By
   definition, every $\rh$-submodule of $M$ is $\rh$-Noetherian,
   and in particular so is ~$\rd\cdot m$. Hence $\rd/\fp$ is
   $\rh$-Noetherian for every $\fp\in\asc(M)$.

   Conversely, suppose that $\rd/\fp$ is $\rh$-Noetherian for
   every $\fp\in\asc(M)$. Since $M$ is $\rd$-Noetherian, there is
   a chain of $\rd$-submodules
   \begin{displaymath}
      0=M_0\subset M_1 \subset \dots \subset M_{r-1}
      \subset M_r=M
   \end{displaymath}
   with $M_j/M_{j-1}\cong \rd/\fq_j$ for $1\le j\le r$, where
   each $\fq_j$ is a prime ideal in $\rd$ that contains some
   $\fp_j\in\asc(M)$ \cite[Corollary 6.2]{Sch}. The surjection
   $\rd/\fp_j\to\rd/\fq_j$ shows that $\rd/\fq_j$ is again
   $\rh$-Noetherian. Recall the fact from module theory that over
   an arbitrary ring if $P\subset Q$ are modules, then $Q$ is
   Noetherian if and only if both $P$ and $Q/P$ are Noetherian.
   Repeated application of this fact shows successively that
   $M_1$, $M_2$, $\dots$, and finally $M_r=M$ are
   $\rh$-Noetherian.
\end{proof}

We will use twice the following variant of the ``determinant
trick,'' which for clarity we state separately.

\begin{lemma}\label{lem:determinant-trick}
   Let $\mathscr{R}$ be a commutative ring with unit and
   $\mathscr{M}$ be an $\mathscr{R}$-module generated by $b_1$,
   $\dots$, $b_k$. Suppose that $r_{ij}\in\mathscr{R}$ for $1\le
   i,j\le k$ are such that $\sum_{j=1}^k r_{ij}b_j=0$ for $1\le
   i\le k$. Then $\det[r_{ij}]$ annihilates $\mathscr{M}$.
\end{lemma}

\begin{proof}
   Consider the matrix $A=[r_{ij}]$ acting on column vectors in
   ~$\mathscr{M}^k$. If $\mathbf{b}\in\mathscr{M}^k$ has $j$th
   entry $b_j$, then $A\mathbf{b}=0$. Multiplying this equation
   by the adjugate matrix of $A$ gives $(\det A)\mathbf{b}=0$.
   Hence $(\det A)b_j=0$ for $1\le j\le k$, so that $\det A$
   annihilates ~$\mathscr{M}$.
\end{proof}

In order to show that expansiveness of $\al_M$ along $H$ implies
that $M$ is $\rh$-Noetherian, we need one further algebraic
result. Recall that modules over Noetherian rings are Noetherian
exactly when they are finitely generated. This fails for
non-Noetherian rings. For example, if $\rh$ is non-Noetherian,
then $\rh$ is finitely generated over itself (by 1), yet the
ideal $\ZZ[u^{\bn}:\bn\cdot\vH<0]$ is not finitely generated over
~$\rh$. Nevertheless, in our situation there is a reasonable
substitute.

\begin{lemma}\label{lem:finitely-generated}
   Let $M$ be an $\rd$-module and $H\in\Hd$. Then $M$ is
   $\rh$-Noetherian if and only if $M$ is finitely generated over
   $\rh$.
\end{lemma}

\begin{proof}
   Noetherian modules over arbitrary rings are always finitely
   generated.

   Conversely, suppose that $M$ is finitely generated over $\rh$,
   say by $m_1$, $\dots$, $m_r$. Fix $\bk\in\zd\smallsetminus H$. Since
   $M$ is an $\rd$-module, we can find $f_{ij}(u)\in\rh$ such
   that
   \begin{displaymath}
      u^{\bk}\cdot m_j=\sum_{i=1}^r f_{ij}(u)\cdot m_i.
   \end{displaymath}
   If $I$ denotes the $r\times r$ identity matrix and
   $F=[f_{ij}(u)]$, then Lemma~ \ref{lem:determinant-trick}
   shows that $\det(u^{\bk}I-F)$ annihilates $M$.
   Multiplying this determinant by $u^{-(r-1)\bk}$ shows that
   there is an element in $\rd$ of the form $u^{\bk}-f(u)$, with
   $f(u)\in\rh$, that annihilates~$M$.

   Suppose that $N$ is an $\rh$-submodule of ~$M$. For every $n\in
   N\subset M$ we have that $(u^{\bk}-f(u))\cdot n=0$, so that
   $u^{\bk}\cdot n=f(u)\cdot n\in N$. Hence $N$ is closed under
   the subring of $\rd$ generated by $\rh$ and ~$u^{\bk}$, which
   is all of $\rd$. Thus every $\rh$-submodule of $M$ is also an
   $\rd$-submodule. Since $M$ is finitely generated over $\rh$,
   it is finitely generated over ~$\rd$, hence $M$ is
   $\rd$-Noetherian. It then follows
   that $M$ is also $\rh$-Noetherian (since every $\rh$-submodule
   is also an $\rd$-submodule).
\end{proof}

The proof of Lemma \ref{lem:finitely-generated} also establishes
the following useful result.

\begin{lemma}\label{lem:noetherian-polynomial}
   Let $M$ be a Noetherian $\rd$-module, $H\in\Hd$, and
   $\bk\in\zd\smallsetminus H$. Then $M$ is $R_H$-Noetherian if
   and only if there is a polynomial of the form $u^{\bk}-f(u)$
   with $f(u)\in R_H$ that annihilates $M$. 
\end{lemma}

\begin{remark}\label{rem:Nn-is-closed}
   If $M$ is a Noetherian $\rd$-module, Lemma
   \ref{lem:noetherian-polynomial} shows that
   $\{H\in\Hd:\text{$M$ is $R_H$-Noetherian}\}$ is an open subset
   of $\Hd$.
   To see this, suppose that $M$ is $R_H$-Noetherian
   for some $H\in\Hd$. Fix $\bk\in\zd\smallsetminus H$. Applying
   Lemma \ref{lem:noetherian-polynomial} with $\bk$ replaced by
   $2\bk$, we see there is an $f(u)\in R_H$ such that
   $u^{2\bk}-f(u)$ annihilates $M$. Then $u^{\bk}-u^{-\bk}f(u)$
   also annihilates $M$. For $H'$ sufficiently close to $H$ it
   follows that $u^{-\bk}f(u)$ is also in $R_{H'}$ and that
   $\bk\in\zd\smallsetminus H'$, and hence by another application
   of Lemma \ref{lem:noetherian-polynomial} we find that $M$ is
   $R_{H'}$-Noetherian as well.
\end{remark}

\begin{lemma}\label{lem:exp-implies-noetherian}
   Let $M$ be an $\rd$-module and $H\in\Hd$. If $\al_M$
   is expansive along $H$ then $M$ is $\rh$-Noetherian.
\end{lemma}

\begin{proof}
   If $M$ is not $\rh$-Noetherian, then by Lemma~
   \ref{lem:finitely-generated} it is not finitely generated over
   ~$\rh$. Hence there is a strictly increasing sequence of proper
   $\rh$-submodules $M_1\subsetneqq M_2\subsetneqq\dots$ with
   $\bigcup_{j=1}^\infty M_j=M$. Let $X_j^{}=M_j^\perp\subset
   X_M$. Then each $X_j$ is compact, $X_1\supsetneqq
   X_2\supsetneqq\dots$, and $\bigcap_{j=1}^\infty X_j=\{0_X\}$.
   Hence $\diam(X_j)\to0$. Furthermore, since $M_j$ is an
   $\rh$-module, $\al_M^{\bn}(X_j)\subset X_j$ for every
   $\bn\in\hz$. 
   
   Let $\varepsilon>0$. Choose $j_0$ such that
   $\diam(X_{j_0})<\varepsilon$, and pick $0\ne x\in X_{j_0}$. Then
   $\al_M^{\bn}(x)\in X_{j_0}$ for all $\bn\in\hz$, so that
   $\sup\{\rho(\al_M^{\bn}(x),0_X):\bn\in\hz\}<\varepsilon$.
   Since $\varepsilon$ was arbitrary, we see that $\al_M$ is not
   expansive along ~$H$. This contradiction proves that $M$ is
   $\rh$-Noetherian.
\end{proof}

Next we turn to the one-sided version of the variety condition.
Define $\log|\bz|$ for $\bz=(z_1,\dots,z_d)\in(\CC^\times)^d$ by
\begin{displaymath}
   \log |\bz| = (\log|z_1|,\dots,\log|z_d|)\in\RR^d.
\end{displaymath}
If $\fa$ is an ideal in $\rd$, put
\begin{displaymath}
   \log|\V(\fa)|=\{\log|\bz|:\bz\in\V(\fa)\}\subset\RR^d.
\end{displaymath}
When $\fa=\< f\>$ is principal, the set
$\log|V(\fa)|=\log|\V(f)|$ was investigated in \cite[\S6.1]{GKZ},
where it is called the \textit{amoeba} of $f$ (turn to Figure
\ref{fig:two-variable-variety} to see why). For example, they
show that the connected components of the complement of
$\log|\V(f)|$ are all convex sets that are in one-to-one
correspondence with the distinct domains of convergence of
Laurent expansions of $1/f$.

For $\bv\in\sphd$, let $[0,\infty)\bv=\{t\bv:t\ge0\}$ denote the
ray in $\RR^d$ through ~$\bv$.

\begin{proposition}\label{prop:principal-case}
   Let $\fa$ be an ideal in $\rd$ and $H\in\Hd$ with outward
   unit normal~ $\vH$. Then $\arda$ is expansive along $H$ if and
   only if $\rda$ is $\rh$-Noetherian and
   $[0,\infty)\vH\cap\log|\V(\fa)|=\emptyset$. 
\end{proposition}

\begin{proof}
   Let $\sigma$ denote the $\zd$-shift action on ~$\TT^{\zd}$,
   where $(\sigma^{\bn}x)_{\bm}=x_{\bn+\bm}$. For
   $f=\sum_{\bn\in\zd} c_f(\bn)u^{\bn}\in\rd$ put
   $f(\sigma)=\sum_{\bn\in\zd}c_f(\bn)\sigma^{\bn}$. Then
   $X_{\rda}$ is the closed shift-invariant subgroup of
   $\TT^{\zd}$ given by
   \begin{displaymath}
     X_{\rda}=\{x\in\TT^{\zd}:f(\sigma)x=0\text{ for all
       $f\in\fa$}\}.
   \end{displaymath}
   If $f+\fa\in\rd/\fa$ and $x\in X_{\rda}$, the duality pairing
   is given by the formula
   \begin{displaymath}
      \< f,x \> =\exp\bigl[ 2\pi i\bigl(f(\sigma)x\bigr)
      _{\mathbf{0}}\bigr].
   \end{displaymath}
   
   First suppose that $\arda$ is expansive along ~$H$. Lemma~
   \ref{lem:exp-implies-noetherian} shows that $\rda$ is
   $\rh$-Noetherian. If
   $[0,\infty)\vH\cap\log|\V(\fa)|\ne\emptyset$, choose
   $\bz\in\V(\fa)$ such that $(\log|\bz|)\cdot\bn\le0$ for all
   $\bn\in\hz$. Consider $\CC^{\zd}$ together with the
   $\zd$-shift action~$\sigma$, and the point $w\in\CC^{\zd}$
   defined by $w_{\bn}=\bz^{\bn}=z_1^{n_1}\cdots z_d^{n_d}$. For
   every $f\in\rd$ we have $f(\sigma)(w)=f(\bz)w$, and so
   $f(\sigma)(w)=0$ for all $f\in\fa$. Fix $\varepsilon>0$. Then
   $|\varepsilon w_{\bn}|=|\varepsilon\bz^{\bn}|\le\varepsilon$
   for every $\bn\in\hz$. Define $x\in\TT^{\zd}$ by
   $x_{\bn}=\Re(\varepsilon w_{\bn})\pmod1$. Then clearly $x\ne0$, and
   $|x_{\bn}|\le\varepsilon$ for all $\bn\in\hz$. Since
   $\varepsilon$ was arbitrary, this contradicts expansiveness of
   $\arda$ along ~$H$. Hence
   $[0,\infty)\vH\cap\log|\V(\fa)|=\emptyset$.

   Conversely, suppose that $\rda$ is $\rh$-Noetherian. 
   Lemma \ref{lem:noetherian-polynomial} shows that there is a
   polynomial $g$ in $\fa$ of the form $g(u)=1-\sum_{\bn\in
   G}c_g(\bn)u^{\bn}$, where $\bn\cdot\vH<0$ for all $\bn\in G$.
   Now $\fa$ is finitely generated over $\rd$, say by $f_1$,
   $\dots$, $f_r$.

   For $h=\sum_{\bn\in\zd}c_h(\bn)u^{\bn}\in\rd$ put
   $\|h\|=\sum_{\bn\in\zd}|c_h(\bn)|$. Define
   $\varepsilon=(10\|g\|+10\sum_{j=1}^r \|f_j\|)^{-1}$.
   For $t\in\TT$ let $|t|=\min\{|t-n|:n\in\ZZ\}$.
   We will
   prove that if $x\in X_{\rda}$ and $|x_{\bn}|<\varepsilon$ for
   all $\bn\in\hz$, then $x=0$, which will show that $\arda$ is
   expansive.

   There are two cases to consider: (1) $x_{\bn}=0$ for all
   $\bn\in\hz$, and (2) $x_{\bn}\ne0$ for some $\bn\in\hz$.

   In case (1), choose $\theta>0$ so that $\bn\cdot\vH\le-\theta$
   for all $\bn\in G$. If $\bk\in\zd$ with
   $0<\bk\cdot\vH\le\theta$ then $(\bn+\bk)\cdot\vH\le0$ for
   $\bn\in G$. Now $u^{\bk}g(u)=u^{\bk}-\sum_{\bn\in
   G}c_g(\bn)u^{\bn+\bk}\in\fa$, so that $x_{\bk}=\sum_{\bn\in
   G}c_g(\bn)x_{\bn+\bk}=0$ since $x_{\bn+\bk}=0$ for all $\bn\in
   G$ by assumption. Hence $x_{\bk}=0$ for all $\bk$ with
   $0<\bk\cdot\vH\le\theta$. Repeating this argument shows that
   $x_{\bk}=0$ if $\bk\cdot\vH\le2\theta$, then if
   $\bk\cdot\vH\le3\theta$, and so on, which proves that $x=0$.
   
   We now turn to case (2). We will show that there is a point in
   $[0,\infty)\vH\cap\log|\V(\fa)|$. Consider the Banach space
   $\ell^\infty(\hz)$ of all bounded complex-valued functions on
   ~$\hz$.  For every $\bn\in\hz$ the operator $U_{\bn}$ defined
   by $(U_{\bn}w)_{\bm}=w_{\bm+\bn}$ maps $\ell^\infty(\hz)$ to
   itself, and clearly $\|U_{\bn}\|\le1$. For
   $f=\sum_{\bn\in\hz}c_f(\bn)u^{\bn}\in\rh$ define the operator
   $f(U)=\sum_{\bn\in\hz}c_f(\bn)U_{\bn}$. Put
   \begin{align*}
      W=\{\,w\in\ell^\infty(\hz):{} U_{\bm}&f_j(U)w=0 \text{ and }
          U_{\bm}g(U)w=0 \\
          &\text{ for $1\le j\le r$ and all $\bm\in\hz$}\}.
   \end{align*}
   Clearly $W$ is closed and mapped to itself by $U_{\bn}$ for
   every $\bn\in\hz$. We claim that $W$ is nontrivial. For define
   $w\in W$ by taking $w_{\bn}$ to be the unique number in
   $(-\varepsilon,\varepsilon)$ for which $x_{\bn}\equiv
   w_{\bn}\pmod1$. Since
   \begin{displaymath}
      1=\< u^{\bm}f_j(u),x\>=
      \exp\bigl[ 2\pi i \bigl( \sigma^{\bm}f_j(\sigma)x\bigr)_{\mathbf{0}}
         \bigr],
   \end{displaymath}
   if follows that
   $\bigl(\sigma^{\bm}f_j(\sigma)w\bigr)_{\mathbf{0}}=\sum
   c_{f_j}(\bn)w_{\bm+\bn}\in\ZZ$ for all $\bm\in\hz$. Our size
   condition on $x$ coupled with $u^{\bm}f_j(u)\in\rh$ shows that
   $\sum c_{f_j}(\bn)w_{\bm+\bn}=0$ for $1\le j\le r$ and
   $\bm\in\hz$. The same argument works for $g$, proving that $w$
   is a nonzero element of~$W$.

   For each $\bn\in\hz$ let $V_{\bn}$ be the restriction of
   $U_{\bn}$ to~ $W$. Let $\SA$ be the commutative Banach algebra
   of bounded operators on $W$ generated by
   $\{V_{\bn}:\bn\in\hz\}$. The theory of commutative Banach
   algebras shows that there is a (nonzero) complex
   homomorphism $\omega\colon\SA\to\CC$, and that
   $|\omega(V)|\le\|V\|$ for all $V\in\SA$. Let
   $a_{\bn}=\omega(V_{\bn})$ for all $\bn\in\hz$. Then
   $\bn\mapsto a_{\bn}$ is a homomorphism from the monoid $\hz$
   to the multiplicative monoid $\CC$. It follows that there is 
   $\bz=(z_1,\dots,z_d)\in\CC^d$ such that $a_{\bn}=\bz^{\bn}$
   for all $\bn\in\hz$. Also, $f_j(V)=0$ on $W$, so
   applying $\omega$ shows that $f_j(\bz)=g(\bz)=0$.

   We claim that $z_j\ne0$ for $1\le j\le d$. Since $\omega\ne0$,
   there is a $\bk\in\hz$ such that $a_{\bk}\ne0$. Suppose that
   $a_{\bn}=0$ for all $\bn\in\hz$ with $\bn\cdot\vH<0$. If
   $\bm\in\partial H\cap\zd$, then $a_{\bm}=\sum_{\bn\in
   G}c_g(\bn)a_{\bm+\bn}=0$, so that $a$ is also zero on $\hz$, a
   contradiction. Hence there is an $\bn\in\hz$ with
   $\bn\cdot\vH<0$ and $a_{\bn}\ne0$. It then follows from
   multiplicativity of the $a_{\bn}$ that $a_{\bn}\ne0$ for all
   $\bn\in\hz$. Hence $z_j\ne0$ for all $j$, and so
   $\bz\in\V(\fa)$. Finally, for all $\bn\in\hz$ we have that
   \begin{displaymath}
      |\bz^{\bn}|=|a_{\bn}|=|\omega(V_{\bn})|\le\|V_{\bn}\|\le1,
   \end{displaymath}
   so that $\bn\cdot\log|\bz|\le0$, proving that
   $\log|\bz|\in[0,\infty)\vH$. 
\end{proof}

We will need the following fact when dealing with general
$\rd$-modules. 

\begin{lemma}\label{lem:build-up}
   Let $\al$ be an algebraic $\zd$-action on $X$,  $Y$ be an
   $\al$-invariant compact subgroup of $X$, and $H\in\Hd$. If
   $\al_Y$ and $\al_{X/Y}$ are both expansive along $H$, then
   so is ~$\al$.
\end{lemma}

\begin{proof}
   The hypothesis shows that there is a neighborhood $U$ of $0_X$ 
   in $X$ such that
   $\bigcap_{\bn\in\hz}\al^{\bn}_{X/Y}(U+Y)=\{0_{X/Y}\}$ and
   $\bigcap_{\bn\in\hz} \al_Y^{\bn}(U\cap Y)=\{0_X\}$. These
   imply that $\bigcap_{\bn\in\hz}\al^{\bn}(U)=\{0_X\}$, so that
   $\al$ is expansive.
\end{proof}

We are now ready to characterize the expansive half-spaces for
algebraic $\zd$-actions. First observe that by Theorem~
\ref{thm:expansive} if an $\rd$-module $M$ is not
$\rd$-Noetherian, then $\al_M$ is not expansive, and so all
subsets of $\RR^d$ are also nonexpansive and we are done. Thus we
can assume the modules considered are Noetherian over ~$\rd$.

\begin{theorem}\label{thm:characterization}
   Let $M$ be a Noetherian $\rd$-module, $\al_M$ be the
   corresponding algebraic $\zd$-action, and $H\in\Hd$. Then the
   following are equivalent.
   \begin{enumerate}
     \item[(1)] $\al_M$ is expansive along $H$.
     \item[(2)] $\ardp$ is expansive along $H$ for every
      $\fp\in\asc(M)$. 
     \item[(3)] $\rdp$ is $\rh$-Noetherian and
      $[0,\infty)\vH\cap\log|\V(\fp)|=\emptyset$ for every
      $\fp\in\asc(M)$. 
   \end{enumerate}
\end{theorem}

\begin{proof}
   Proposition~ \ref{prop:principal-case} shows that $(2)\ifff(3)$.

   $(3)\implies(1)$: There is a chain of $\rd$-submodules
   $0=M_0\subset M_1\subset\dots\subset M_r=M$ such that
   $M_j/M_{j-1}\cong \rd/\fq_j$, where $\fq_j$ is a prime ideal
   in $\rd$ containing some $\fp_j\in\asc(M)$. Since
   $X_{\rd/\fq_j}\subset X_{\rd/\fp_j}$, we see that
   $\al_{\rd/\fq_j}$ is expansive along ~$H$. Put
   $X_j^{}=M_j^\perp\subset X_M$. Then $X_M=X_0\supset X_1\supset\dots
   X_r=0$, and $X_{j-1}/X_j\cong X_{\rd/\fq_j}$. Repeated
   application of Lemma~ \ref{lem:build-up} shows successively
   that $\al_{X_{r-1}}$, $\al_{X_{r-2}}$, $\dots$,
   $\al_{X_0}=\al_M$ are all expansive along ~$H$.

   $(1)\implies(3)$: Suppose that $\al_M$ is expansive along ~$H$.
   By Lemma~ \ref{lem:exp-implies-noetherian}, $M$ is
   $\rh$-Noetherian, so that by Lemma~ \ref{lem:Rh-noetherian} we
   have that $\rdp$ is $\rh$-Noetherian for every $\fp\in\asc(M)$.

   Suppose that there is a $\fp\in\asc(M)$ for which
   $[0,\infty)\vH\cap \log|\V(\fp)|\ne\emptyset$. Fix 
   $\bz\in\V(\fp)$ with $\log|\bz|=t\vH$ for some $t\ge0$. 
   Consider $\CC$ as an $\rd$-module using the action $\theta$
   defined by $\theta^{\bn}(\xi)={\bz}^{\bn}\xi$ for all
   $\xi\in\CC$. We will construct an $\rd$-homomorphism
   $\psi\colon M\to\CC$.
   
   To construct $\psi$, first note that $M$ is $\rd$-Noetherian
   by assumption, hence finitely generated over
   ~$\rd$. Choose generators $m_1$, $m_2$, $\dots$, $m_k$, and
   define the surjective map $\zeta\colon R_d^k\to M$ by
   $\zeta(f_1,\dots,f_k)=f_1\cdot m_1+\dots+f_k\cdot m_k$. Let
   $K=\kernel \zeta$. Define $\phi\colon R_d^k\to \CC^k$ by
   $\phi(f_1,\dots,f_k)=(f_1(\bz),\dots,f_k(\bz))$. We claim that
   the dimension of the complex vector space generated by
   $\phi(K)$ is strictly less than ~$k$. For if not, there are
   elements $\bf^{(1)}=(f_1^{(1)},\dots,f_k^{(1)})$, $\dots$,
   $\bf^{(k)}=(f_1^{(k)},\dots,f_k^{(k)})$ in $K$ such that
   $\phi(\bf^{(1)})$, $\dots$, $\phi(\bf^{(k)})$ are linearly
   independent in $\CC^k$. Denote the $k\times k$ matrix
   $[f_j^{(i)}]$ by $F$. Since each $\bf^{(i)}\in K$, Lemma~
   \ref{lem:determinant-trick} shows that $\det F$ annihilates
   ~$M$, and in particular $\det F\in\fp$. Thus
   \begin{displaymath}
      0=(\det F)(\bz)=\det
      [f^{(i)}_j(\bz)]=\det[\phi(\bf^{(1)}),\dots, \phi(\bf^{(k)})],
   \end{displaymath}
   contradicting linear independence of the ~$\phi(\bf^{(i)})$.
   This proves that $\phi(K)$ generates a proper complex vector
   subspace ~$L$ of ~$\CC^k$. Since $\phi(R_d^k)$ generates all of
   $\CC^k$, it follows that the corresponding quotient map
   $\widetilde{\phi}\colon M\cong R_d^k/K\to\CC^k/L$ is nonzero.
   We can therefore compose $\widetilde{\phi}$ with a projection
   of $\CC^k/L$ to a 1-dimensional subspace so that the
   composition is still a nonzero $\rd$-homomorphism. The result
   is the desired $\psi\colon M\to\CC$.

   By construction, $\psi(f_j \cdot m_j)=c_j f_j(\bz)$. After
   multiplying (if necessary) by a constant we can assume that
   $|c_j|<\varepsilon$ and that the point $x\in X_M=\widehat{M}$
   defined by $x(u^{\bn}\cdot m_j)=\exp\bigl[2\pi i
   \Re(c_j\,\bz^{\bn})\bigr]$ is not trivial. Since
   $\al^{\bn}(x)$ is close to $0_X$ for all $\bn\in\hz$,
   the half-space $H$ is not expansive.
\end{proof}

\begin{remark}\label{rem:recover-original}
   Roughly speaking, we can recover the statement (and proof) of
   Theorem \ref{thm:expansive} from the preceding theorem by
   using $\bv=\vH=0$. For then $[0,\infty)\vH=\{\mathbf{0}\}$, so
   that $[0,\infty)\vH\cap \log|\V(\fp)|=\emptyset$ if and only
   if $\V(\fp)\cap\SS^d=\emptyset$, and
   $H_{\bv}=\{\bx\in\RR^d:\bx\cdot\bv\le0\}=\RR^d$. 
\end{remark}

The content of Theorem~ \ref{thm:characterization} can be summarized by 
the equalities
\begin{displaymath}
   \E(\al_M)=\bigcap_{\fp\in\asc(M)}\E(\ardp),\quad
   \N(\al_M)=\bigcup_{\fp\in\asc(M)}\N(\ardp).
\end{displaymath}

Our arguments would be substantially simplified if we knew that
the converse of Lemma~ \ref{lem:build-up} were true, namely that
the quotient of an algebraic $\zd$-action that is expansive along
~$H$ is also expansive along ~$H$. Although this is correct, the
only proof we know makes full use of Theorem~
\ref{thm:characterization}. 

\begin{proposition}\label{prop:quotients}
   Let $\al$ be an algebraic $\zd$-action on $X$, $Y$ be a
   compact $\al$-invariant subgroup of $X$, and $H\in\Hd$. Then
   $\al$ is expansive along ~$H$ if and only if both $\al_Y$ and
   $\al_{X/Y}$ are expansive along ~$H$.
\end{proposition}

\begin{proof}
   Lemma~ \ref{lem:build-up} proves one direction. For the other,
   suppose $\al$ is expansive along ~$H$. Then trivially $\al_Y$
   is expansive along ~$H$. Let $N=Y^\perp\subset M$, so that
   $X/Y=X_N$. Now $\asc(N)\subset\asc(M)$, so Theorem~
   \ref{thm:characterization} shows that $\al_N=\al_{X/Y}$ is
   expansive along ~$H$.
\end{proof}

\section{Examples}\label{sec:examples}

This section contains a number of examples to illustrate the
above ideas. One noteworthy feature is the elegant way in which
the non-Noetherian and variety pieces of the nonexpansive set fit
together. 

Our examples involve three or fewer variables, so for notational
simplicity we use $u$, $v$, $w$, instead of $u_1$, $u_2$, $u_3$.
Using the correspondence $\Hd \leftrightarrow \sphd$ given by
$H\leftrightarrow \vH$, we identify subsets of $\Hd$ with the
corresponding subsets of $\sphd$ for ease of visualization.
Using this convention, for an ideal $\fa\in\rd$ we put
\begin{align*}
   \Nn(\arda)&=\{\bv\in\sphd:\text{$\rd/\fa$ is not
     $R_{H_{\bv}}$-Noetherian}\}, \\
   \Nv(\arda)&=\{\bv\in\sphd:[0,\infty)\bv\cap
     \log|\V(\fa)|\ne\emptyset\}.
\end{align*}
Observe that $\Nv(\arda)$ is the radial projection of $\log
|\V(\fa)|$ to $\sphd$.  By Proposition~ \ref{prop:principal-case},
\begin{displaymath}
   \N(\arda)=\Nn(\arda)\cup\Nv(\arda).
\end{displaymath}
As before, in the case of a principal ideal $\< f\>$ in $\rd$ we abbreviate
$\V(\< f\>)$ to ~$\V(f)$.

\begin{example}\label{exam:one-variable-principal}
   (\textit{One variable, principal ideal}.)  Let $0\ne f(u)\in
   R_1$, which we assume is not a unit in $R_1$ (i.e., not
   $\pm1$ times a monomial). Multiplying
   $f(u)$ by a monomial if necessary, we can also assume that
   $f(u)=c_ru^r+c_{r-1}u^{r-1}+ \dots+ c_1u+c_0$, where
   $c_j\in\ZZ$ and $c_rc_0\ne0$. The unit ``sphere'' in $\RR$ is
   $\sph_0=\{1,-1\}$. By Lemma \ref{lem:noetherian-polynomial},
   $1\in\Nn(\arof)$ if and only if $|c_r|>1$. If $|c_r|=1$, then
   $f(u)=\pm\prod_{j=1}^r(u-\lambda_j)$, so that
   $\prod_{j=1}^r|\lambda_j|=|c_0|\ge1$.
   Hence $a=\max_{1\le j\le r}|\lambda_j|\ge1$. If $a=1$, then
   $\arof$ is not expansive, 
   and so trivially
   $1\in\N(\arof)$. If $a>1$, then
   $[0,\infty)\cdot1\cap \log|\V(f)|\ne\emptyset$, and so
   $1\in\Nv(\arof)$. In all cases we conclude that
   $1\in\N(\arof)$. Similarly, $-1\in\N(\arof)$. Thus
   $\N(\arof)=\sph_0$.

   The following polynomials illustrate some possible
   combinations of $\Nn$ and $\Nv$.
   \begin{enumerate}
     \item[(a)] $f(u)=u^2-u-1$, $\Nn=\emptyset$ and
      $\Nv=\{1,-1\}$.
     \item[(b)] $f(u)=u-2$, $\Nn=\{-1\}$ and $\Nv=\{1\}$.
     \item[(c)] $f(u)=2$, $\Nn=\{1,-1\}$ and $\Nv=\emptyset$.
     \item[(d)] $f(u)=2u^2-6u+3$, $\Nn=\{1,-1\}$ and $\Nv=\{1,-1\}$.
   \end{enumerate}
\end{example}

\begin{remark}\label{rem:no-one-sided-expansive}
   A result going back to the Ph.D. thesis of Schwartzman shows
   that there are no ``one-sided expansive homeomorphisms''
   except on finite spaces (see \cite[Theorem~ 3.9]{BL} for a
   discussion). From this it follows that if $\alpha$ is an
   algebraic $\ZZ$-action on an infinite group, then
   $\N(\al)=\sph_0$, providing an alternative approach to
   Example~ \ref{exam:one-variable-principal}.
\end{remark}

\begin{example}\label{exam:two-variables-principal}
   (\textit{Two variables, principal ideal}.) Let $f=\sum
   c_f(\bn)u^{\bn}\in R_2$. The Newton polyhedron $\Newt(f)$ of
   $f$ is the convex hull of $\{\bn\in\ZZ^2:c_f(\bn)\ne0\}$. If
   $\Newt(f)$ is a point or line segment, then we are essentially
   reduced to Example~ ~\ref{exam:one-variable-principal}, so we
   assume here that $\Newt(f)$ is $2$-dimensional.
   
   List the vertices of $\Newt(f)$ as $\bn_1$, $\bn_2$, $\dots$,
   $\bn_r$, so that the line segment $[\bn_j,\bn_{j+1}]$ is an
   edge of $\Newt(f)$ (with the convention that
   $\bn_{r+1}=\bn_1$). Let $\bv_j\in\sph_1$ denote the outward
   unit normal vector to $[\bn_j,\bn_{j+1}]$. If $A_j$ denotes
   the open arc from $\bv_{j-1}$ to ~$\bv_{j}$, then $\sph_1$ is
   subdivided into the points ~$\bv_j$ and arcs ~$A_j$ (see
   Figure~ \ref{fig:two-variable-newton}).  This subdivision of
   $\sph_1$ represents the ``spherical dual polygon'' to
   ~$\Newt(f)$, with vertices ~$\bn_j$ of $\Newt(f)$
   corresponding to edges ~$A_j$, and edges $[\bn_j,\bn_{j+1}]$
   in $\Newt(f)$ corresponding to vertices ~$\bv_j$.

   \begin{figure}[htbp]
      \begin{center}
         \scalebox{.9}{\includegraphics{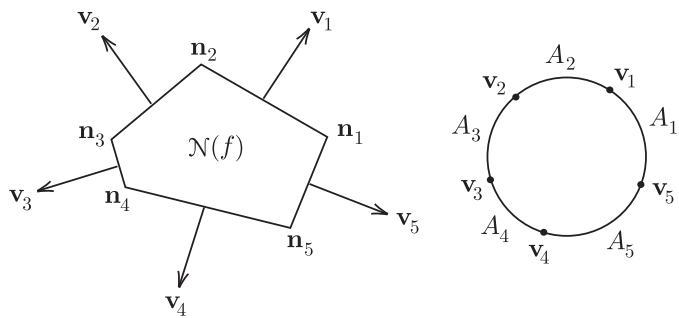}}
         \caption{Newton polygon and its spherical dual}
         \label{fig:two-variable-newton}
      \end{center}
   \end{figure}
   
   By Lemma \ref{lem:noetherian-polynomial}, each
   $\bv_j\in\Nn(\artf)$. If $|c_f(\bn_j)|>1$, we also have that
   $A_j\subset\Nn(\artf)$. If $|c_f(\bn_j)|=1$, then an argument
   using $\V(f)$ similar to that in Example~
   \ref{exam:one-variable-principal} shows that
   $A_j\subset\Nv(\artf)$. Hence in every case
   $A_j\subset\N(\artf)$, so that $\N(\artf)=\sph_1$ (see
   \cite{Scott} for a detailed treatment of this situation).

   We examine two specific polynomials. To describe sets in
   $\sph_1$ we use the notation $\bv_{\theta}=(\cos
   \theta,\sin\theta)$.
   \begin{enumerate}
     \item[(a)] Let $f(u,v)=3+u+v$. Then $\Nn(\artf)=\{\bv_{\theta}:
      \pi\le\theta\le3\pi/2\}\cup\{\bv_{\pi/4}\}$. The set
      $\log|\V(f)|$ is depicted in Figure
      \ref{fig:two-variable-variety}(a), where the boundary
      curves are parameterized by $(\log r,\log|3\pm r|)$ for
      $0<r<\infty$. Projecting this set radially to
      $\sph_1$ shows that
      $\Nv(\artf)=\{\bv_{\theta}:-\pi/2<\theta<\pi\}$, and so
      $\N(\artf)=\sph_1$. 
     \item[(b)] Let $f(u,v)=5+u+u^{-1}+v+v^{-1}$. Then
      $\Nn(\artf)$ consists of just the four points on
      $\sph_1$ corresponding to the four outward unit
      normals of $\Newt(f)$. However, $\log|\V(f)|$, shown in
      Figure \ref{fig:two-variable-variety}(b), covers all
      directions, so that $\Nv(\artf)=\sph_1$.
   \end{enumerate}

   \begin{figure}[htbp]
      \begin{center}
         \scalebox{1}{\includegraphics{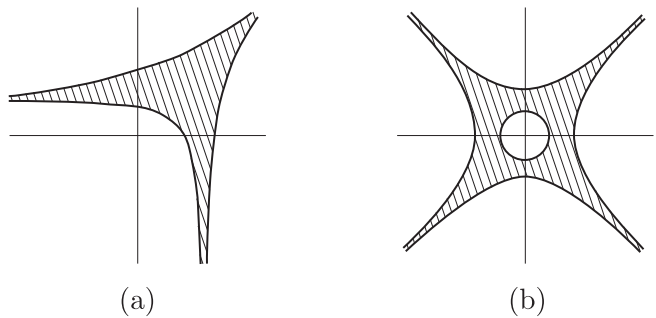}}
         \caption{Logarithmic images of two varieties}
         \label{fig:two-variable-variety}
      \end{center}
   \end{figure}
\end{example}

\begin{example}\label{exam:three-variables-principal}
   (\textit{Three variables, principal ideal}.)
   Let $f=\sum c_f(\bn)u^{\bn}\in R_3$, and $\Newt(f)$ denote the
   Newton polyhedron of ~$f$. If $\Newt(f)$ has dimension $\le2$,
   then we are reduced to either Example~
   \ref{exam:one-variable-principal} or Example~
   \ref{exam:two-variables-principal}, so we will assume that
   $\Newt(f)$ is 3-dimensional.

   We form the ``spherical dual polytope'' to $\Newt(f)$ on ~$\sph_2$
   as follows. For each 2-dimensional face $F$ of $\Newt(f)$
   let $\bv_{\!F}$ be the outward unit normal to ~$F$. If faces $F$
   and $F'$ share a common edge ~$e$, draw a great circle arc
   $A_e$ from $\bv_{\!F}$ to ~$\bv_{\!F'}$. These arcs subdivide
   $\sph_2$ into open regions $B_{\bn}$ corresponding to
   vertices ~$\bn$ of ~$\Newt(f)$ (see Figure~
   \ref{fig:newton-polytope}).

   \begin{figure}[htbp]
      \begin{center}
          \scalebox{1.1}{\includegraphics{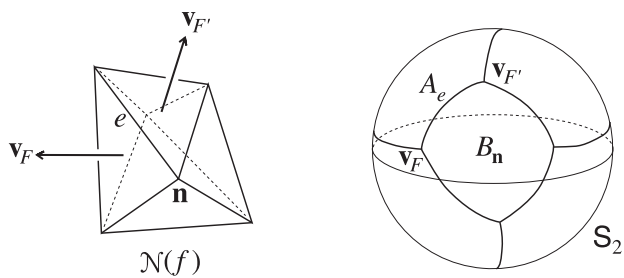}}
         \caption{The Newton polytope and its spherical dual}
         \label{fig:newton-polytope}
      \end{center}
   \end{figure}

   By Lemma \ref{lem:noetherian-polynomial}, each arc
   $A_e\subset\Nn(\arthf)$. If $\bn$ is a vertex of $\Newt(f)$
   with $|c_f(\bn)|>1$, then also $B_{\bn}\subset\Nn(\arthf)$. If
   $|c_f(\bn)|=1$, then $B_{\bn}\cap\Nn(\arthf)=\emptyset$, but
   an analysis similar to the previous examples shows that
   $B_{\bn}\subset\Nv(\arthf)$. Again we find that $\N(\arthf)$
   fills out the entire sphere ~$\sph_2$.
\end{example}

For the principal ideals considered so far we found that the
nonexpansive set is always the whole sphere. This is true in
general.

\begin{proposition}\label{prop:nonexpansive-principal}
   Let $f\in\rd$ generate a proper ideal. Then $\N(\ardf)=\sphd$.
\end{proposition}

\begin{proof}
   We could proceed along the lines of the previous examples,
   which is direct but technically complicated. However, there is
   a short proof using the notion of homoclinic point (which is
   treated in detail in Section~ \ref{sec:homoclinic-group}).
   
   If $\V(f)\cap\SS^d\ne\emptyset$, then $\ardf$ is not expansive
   by Theorem~\ref{thm:expansive}, so that $\N(\ardf)=\sphd$ and
   we are done. So we may suppose that
   $\V(f)\cap\SS^d=\emptyset$. We represent $X_{\rdf}$ as
   $\bigl\{x\in\TT^{\zd}:\sum_{\bn\in\zd}c_f(\bn)x_{\bm+\bn}=0\text{
     for all $\bm\in\zd$}\bigr\}$, so that $\al_{\rdf}$ is the
   shift-action on $X_{\rdf}$. For $t\in\TT$ let
   $|t|=\min\{|t-n|:n\in\ZZ\}$. A point $x\in X_{\rdf}$ is
   \textit{homoclinic} if $|x_{\bn}|\to0$ as $\|\bn\|\to\infty$.
   Consider the function $F$ on $\SS^d$ defined by
   $F(s_1,\dots,s_d)=f(s_1^{-1}, \dots,s_d^{-1})$. Then $F$ does
   not vanish on $\SS^d$ by our assumption on ~$\V(f)$. Let
   $(1/F)\sphat\,(\bn)$ be the Fourier transform of $1/F$ at
   $\bn\in\zd$. Define $\xtri$ by setting $\xtri_{\bn}$ to be the
   reduction mod 1 of $(1/F)\sphat\,(\bn)$. By \cite[Lemma~
   4.5]{LS}, $\xtri\in X_{\rdf}$, and $|\xtri_{\bn}|\to0$ as
   $\|\bn\|\to\infty$ by the Riemann-Lebesgue Lemma. Hence
   $\xtri$ is a homoclinic point.
   
   Let $\bv\in\sphd$ be arbitrary. Let $H=H_{\bv}$ and let
   $\varepsilon>0$. Choose $r$ so that
   $|\xtri_{\bn}|<\varepsilon$ for $\|\bn\|>r$. Pick $\bk\in\zd$
   such that $\dist(\bk,H)>r$, and put
   $y=\al^{\bk}_{\rdf}(\xtri)$. Then $|y_{\bn}|<\varepsilon$ for
   all $\bn\in\hz$. Since $\varepsilon$ was arbitrary,
   $\bv\in\N(\ardf)$ for every $\bv\in\sphd$.
\end{proof}

We remark that a proof of the unoriented version of the previous
result, namely that $\N_{d-1}(\ardf)=\G_{d-1}$, can be given
using entropy. For if $\entropy(\ardf)=0$, then $f$ is a product
of generalized cyclotomic polynomials, and then $\ardf$ is not
expansive. If $\entropy(\ardf)>0$, then the entropy of $\ardf$
along every $(d-1)$-plane is infinite, and so no $(d-1)$-plane
can be expansive.

For nonprincipal prime ideals the nonexpansive set can exhibit
more variety. 

\begin{example}\label{exam:ledrappier-set} 
   (\textit{Ledrappier's example}.) We revisit Example~
   \ref{exam:ledrappier}. Let $\fp=\<
   2,1+u+v\>$, which is easily seen to be a prime ideal in
   ~$R_2$. Then $\V(\fp)=\emptyset$, so that $\Nv(\artp)=\emptyset$.
   On the other hand, the outward normals to the sides of the
   Newton polygon $\Newt(1+u+v)$ are nonexpansive, and by Lemma
   \ref{lem:noetherian-polynomial} these are the only
   nonexpansive vectors. Thus
   $\N(\artp)=\{\bv_{\pi/4},\bv_{\pi},\bv_{3\pi/2}\}$. 
\end{example}

\begin{figure}[htbp]
   \begin{center}
      \scalebox{1.0}{\includegraphics{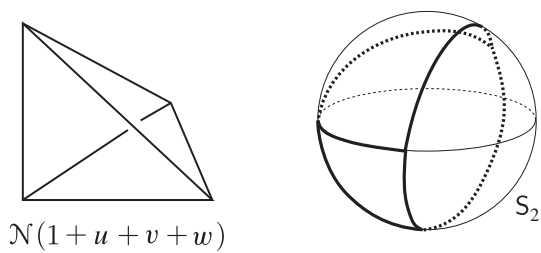}}
      \caption{The nonexpansive set for 3-dimensional Ledrappier}
      \label{fig:3d-ledrappier}
   \end{center}
\end{figure}

\begin{example}\label{exam:3d-ledrappier}
   (\textit{3-dimensional Ledrappier example}.)  Let
   $\fp=\< 2,1+u+v+w\>$, again easily seen to be a
   prime ideal in $R_3$. Then $\V(\fp)=\emptyset$, so that
   $\Nv(\arthp)=\emptyset$. The nonexpansive set was determined
   in \cite[Example~ 2.9]{BL} to be the 1-skeleton of the
   spherical dual to the Newton polytope $\Newt(1+u+v+w)$. This
   is depicted in Figure \ref{fig:3d-ledrappier}, where
   $\Nn(\arthp)$ consists of the six great circle arcs determined
   by the six edges of $\Newt(1+u+v+w)$.
\end{example}

\begin{example}\label{exam:eins-ward}
   (\textit{Nonexpansive set with nonempty interior}.)
   Let $d=3$ and $\fp=\< 1+u+v,w-2\>$. We first prove that
   $\fp$ is a prime ideal in ~$R_3$.

   Define $\phi\colon R_3\to\ztt$ by
   $\phi(f)=f(t,-t-1,2)$. Clearly $\fp\subset\kernel\phi$. Observe
   that $\ztt$ is a subring of $\QQ(t)$, hence an integral
   domain. Define $\psi\colon\ZZ[t]\to R_3/\fp$ by
   $\psi(t)=u+\fp$. Note that
   $\psi(2t(t+1))=2u(u+1)+\fp=-uvw+\fp$, a unit in ~$R_3/\fp$.
   Hence $\phi$ extends uniquely to $\ztt$, and this extension is
   therefore the inverse of the map $R_3/\fp\to\ztt$ induced by
   ~$\phi$. Hence $\fp=\kernel\phi$. Thus $R_3/\fp$ is isomorphic to
   the integral domain $\ztt$, so $\fp$ is prime. (We are
   grateful to Paul Smith for showing us this point of view.)

   Next we determine $\Nv(\arthp)$. Since
   $\V(\fp)=\{(z,-z-1,2):z\in\CC\}$, we see that $\log|\V(\fp)|$
   lies in a plane at height $\log~ 2$ above the origin, and in this
   plane it has the shape shown in Figure \ref{fig:eins-ward}(a),
   where the boundary curves are parameterized by $(\log
   r,\log|r\pm1|)$ for $0<r<\infty$.
   When projected radially to ~$\sph_2$, we obtain the set
   in the upper hemisphere shown in Figure
   \ref{fig:eins-ward}(b), with three cusps on the equator.

   \begin{figure}[htbp]
      \begin{center}
          \scalebox{1.0}{\includegraphics{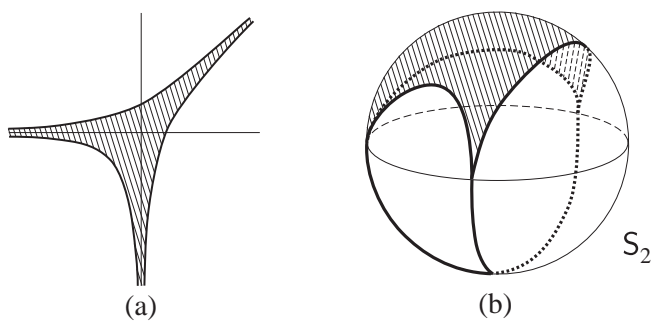}}
         \caption{Nonexpansive set for Example~ \ref{exam:eins-ward}}
         \label{fig:eins-ward}
      \end{center}
   \end{figure}

   Finally, let us compute $\Nn(\arthp)$. Using Lemma
   \ref{lem:noetherian-polynomial}, the polynomial $w-2\in\fp$
   shows that the open upper hemisphere in ~$\sph_2$ is disjoint
   from $\Nn(\arthp)$. Furthermore, $1+u+v\in\fp$ shows that no
   points in the lower hemisphere are in $\Nn(\arthp)$
   either, with the possible exceptions of those on the three quarter
   meridians shown in Figure \ref{fig:eins-ward}(b). We will show
   that each of these quarter meridians is contained in
   $\Nn(\arthp)$, so that they, combined with $\Nv(\arthp)$ in
   the upper hemisphere, comprise all of ~$\N(\arthp)$.
   
   We will treat the meridian from $(0,-1,0)$ to $(0,0,-1)$, the
   other two being similar. Since $\Nn(\arthp)$ is closed by
   Remark ~\ref{rem:Nn-is-closed},
   it is enough to show that unit vectors in the directions
   $(0,-a,-b)$ are non-Noetherian, where $a$ and $b$ are positive
   integers.  Let $H\in\mathsf{H}_3$ be
   $\{\bx\in\RR^3:\bx\cdot(0,a,b)\ge0\}$.  Using the isomorphism
   $\phi\colon R_3/\fp\to\ztt$ above, the subring $R_H$ is mapped
   to $\ZZ[t^{\pm1}, (-t-1)^m2^n:am+bn\ge0]$. Then $R_3/\fp$ is
   Noetherian over $R_H$ if and only if $\ztt$ is finitely
   generated over $\ZZ[t^{\pm1}, (-t-1)^m2^n:am+bn\ge0]$.
   By Lemma \ref{lem:noetherian-polynomial}, this is
   clearly equivalent to whether we can write $1$ as a
   combination, using coefficients in $\ZZ[t^{\pm1}]$, of
   expressions of the form $(-t-1)^m2^n$, where $am+bn>0$.

   Suppose this to be the case, so that
   \begin{equation}\label{eqn:two-adic}
      1=\sum_{(m,n)\in F} f_{mn}(t)(-t-1)^m\,2^n,
   \end{equation}
   where $f_{mn}(t)\in\ZZ[t^{\pm1}]$ and $F$ is a finite set of
   $(m,n)\in\ZZ^2$ for which $am+bn>0$. Let $|\cdot|_2$ denote
   the extension of the $2$-adic norm on $\QQ$ to $\QQ(2^{1/b})$.
   Substitute $t=2^{a/b}-1$ in \eqref{eqn:two-adic}. Since
   $|2^{a/b}-1|_2=1$, it follows that $|f_{mn}(2^{a/b}-1)|_2\le1$.
   Hence
   \begin{align*}
      1=|1|_2 &=\biggl| \sum_{(m,n)\in F} f_{mn}(2^{a/b}-1)
                   (-2^{a/b})^m\, 2^n \biggr|_2 \\
              &\le \max_{(m,n)\in F}\, \bigl| f_{mn}(2^{a/b}-1)\bigr|_2
                   \,\bigl| -2^{a/b}\bigr|_2^m \, \bigl|2\bigr|_2^n \\
              &\le \max_{(m,n)\in F} \,2^{-(am+bn)/b} < 1.
   \end{align*}
   This contradiction shows that \eqref{eqn:two-adic} is
   impossible, so that each rational direction $(0,-a,-b)$ is
   non-Noetherian. See also Example \ref{exam:better-Nn} for
   another approach. 
\end{example}

\begin{remark}\label{rem:bergman}
   In the previous example, and in many others, the
   non-Noetherian and variety parts of the nonexpansive set are
   ``glued together'' along a set which can be described by the
   asymptotic behavior of the logarithmic image of the variety.
   For instance, in Example \ref{exam:eins-ward}, the variety
   part $\Nv(\alpha_{R_3/\fp})$ in the upper hemisphere and the
   non-Noetherian part $\Nn(\alpha_{R_3/\fp})$ in the lower
   hemisphere are glued together at the three cusp points on the
   equator, which are also the three asymptotic directions of
   $\log|\V(\fp)|$. This illustrates a general phenomenon of the
   \textit{logarithmic limit set of an algebraic variety},
   introduced by Bergman \cite{Bergman}. He showed that this set
   is always contained in a finite union of lower-dimensional
   great spheres. We point out that this set has an alternative
   description as the set of half-spaces $H$ for which
   $(\rdp)\otimes\QQ$ is not Noetherian over $R_H\otimes\QQ$, so in
   this sense is also a limiting part of the non-Noetherian set
   (we are grateful to Bernd Sturmfels for pointing this out to us).
\end{remark}

\section{Further analysis of nonexpansive sets}
\label{sec:further-analysis}

Let $M$ be an $\rd$-module, which we will assume throughout this
section to be Noetherian. In Section~ \ref{sec:characterization} we
showed that
\begin{equation}\label{eqn:N(aM)}
   \N(\al_M)=\bigcup_{\fp\in\asc(M)}\N(\ardp),
\end{equation}
where
\begin{equation}\label{eqn:N(aRdp)}
   \N(\ardp)=\Nn(\ardp)\cup\Nv(\ardp).
\end{equation}
In this section we refine our analysis of ~$\N(\al_M)$, showing
that only isolated primes in ~$\asc(M)$ are relevant. We provide
an algorithm using Gr\"obner bases to compute $\Nn(\arda)$ for
ideals $\fa$ in $\rd$, and show how this can be used to compute
$\Nn(\al_M)$ using the Fitting ideal of $M$. The question of when
either $\Nn(\al_M)$ or $\Nv(\al_M)$ is empty is answered in
Proposition \ref{prop:empty}.

We begin by extending the definitions of $\Nn$ and $\Nv$.

\begin{definition}\label{def:N-for-modules}
   For $M$ a Noetherian $\rd$-module, define
   \begin{equation}\label{eqn:extended-defs}
      \Nn(\al_M)=\bigcup_{\fp\in\asc(M)} \Nn(\ardp)
      \text{\quad and\quad}
      \Nv(\al_M)=\bigcup_{\fp\in\asc(M)} \Nv(\ardp).
   \end{equation}
\end{definition}

Hence by \eqref{eqn:N(aM)} and \eqref{eqn:N(aRdp)},
\begin{equation}\label{eqn:N-as-union}
   \N(\al_M)=\Nn(\al_M)\cup\Nv(\al_M).
\end{equation}

\begin{remark}\label{rem:non-Noetherian-M}
   Alternatively, we could have defined $\Nn(\al_M)$ to be
   $\{H\in\Hd: \text{$M$ is not $R_H$-Noetherian}\}$. For $M$ has
   a prime filtration $0=M_0\subset M_1\subset\dots\subset M_r=M$
   with $M_j/M_{j-1}\cong\rd/\fq_j$, where $\fq_j$ is a prime
   ideal containing some $\fp_j\in\asc(M)$, and furthermore every
   $\fp\in\asc(M)$ occurs as some $\fq_j$ \cite[p.\ 93]{Eis}.
   Then $M$ fails to be $R_H$-Noetherian if and only if some
   quotient $M_j/M_{j-1}$ is not $R_H$-Noetherian, establishing
   the equivalence of the two definitions of $\Nn(\al_M)$. Note
   that Remark \ref{rem:Nn-is-closed} shows that $\Nn(\aM)$ is
   always closed.
\end{remark}

Order $\asc(M)$ by inclusion, and let $\min(M)$ denote the subset
of minimal elements of $\asc(M)$. Members of $\min(M)$ are called
the \textit{isolated primes} for $M$. They play an essential role
in the primary decomposition of $M$, and also govern the
expansive subdynamics of $\al_M$

\begin{proposition}\label{prop:min-primes}
   Let $M$ be a Noetherian $\rd$-module. Then
   \begin{equation}\label{eqn:min-primes}
      \Nn(\al_M)=\bigcup_{\fp\in\min(M)} \Nn(\ardp),
      \quad
      \Nv(\al_M)=\bigcup_{\fp\in\min(M)} \Nv(\ardp),
   \end{equation}
   and so
   \begin{displaymath}
      \N(\al_M)=\bigcup_{\fp\in\min(M)} \N(\ardp).
   \end{displaymath}
\end{proposition}

\begin{proof}
   Let $\fq\in\asc(M)$ and choose $\fp\in\min(M)$ with
   $\fp\subset\fq$. The natural surjection $\rd/\fp\to\rd/\fq$
   shows that, for every half-space $H\in\Hd$, if $\rd/\fq$ is
   not $R_H$-Noetherian, then $\rd/\fp$ is not $R_H$-Noetherian.
   Hence $\Nn(\al_{\rd/\fq})\subset\Nn(\ardp)$, establishing the
   first equality in \ref{eqn:min-primes}. For the second
   equality, observe that if $\fp\subset\fq$, then
   $\log|\V(\fq)|\subset \log|\V(\fp)|$.
\end{proof}

By \eqref{eqn:extended-defs} and \eqref{eqn:N-as-union}, we can
compute $\Nv(\aM)$ as the union over all $\fp\in\asc(M)$ of the
radial projections of $\log|\V(\fp)|$ to $\sphd$. We next give
two approaches to computing $\Nn(\aM)$. We first introduce some
convenient terminology.

\begin{definition}\label{def:H-monic}
   Let $f\in\rd$ and $H=\Hv\in\Hd$. The \textit{support of $f$}
   is the set $\{\bn\in\zd:c_f(\bn)\ne0\}$. Say that $f$ has an
   \textit{$H$-exposed vertex $\bn$} if $c_f(\bn)\ne0$ and
   $\bm\cdot\bv<\bn\cdot\bv$ for every
   $\bm\in\supp(f)\smallsetminus\{\bn\}$.  If $f$ has an
   $H$-exposed vertex $\bn$ with $c_f(\bn)=1$, then $f$ is called
   \textit{$H$-monic}.
\end{definition}

Note that Lemma \ref{lem:noetherian-polynomial} shows that $M$ is
$R_H$-Noetherian if and only if there is an $H$-monic polynomial that
annihilates~$M$. In particular, if $\fa$ is an ideal in $\rd$,
then $\rd/\fa$ is $R_H$-Noetherian if and only if $\fa$ contains
an $H$-monic polynomial. Observe that any factor of an $H$-monic
polynomial is $\pm1$ times an $H$-monic polynomial.

We begin with principal ideals $\fa=\<f\>$. By the preceding
paragraph, $\rdf$ is $R_H$-Noetherian if and only if $f$ is $H$-monic. An
obvious extension of terminology in Example~
\ref{exam:three-variables-principal} then shows that $\Nn(\ardf)$
is the $(d-2)$-skeleton of the spherical dual of the Newton
polyhedron $\SN(f)$ of $f$ together with those $(d-1)$-faces of
the spherical dual corresponding to vertices $\bn$ of $\SN(f)$
for which $|c_f(\bn)|>1$.

Now suppose that $\fa=\<f_1,\dots,f_r\>\subset\rd$. For $H\in\Hd$
it is tempting, but wrong, to believe that $\fa$ contains an
$H$-monic polynomial if and only if one of the $f_j$ is $H$-monic.

\begin{example}\label{exam:not-H-monic}
   Let $d=2$, $f=u-2$, $g=v-3$, and $\fa=\<f,g\>$. Take
   $\bv=(-1/\sqrt2,-1/\sqrt2)\in\sph_1$ and $H=\Hv$. Then neither
   $f$ nor $g$ is $H$-monic. However, $f-g=u-v+1\in\fa$ is
   $H$-monic.

   Here
   $\Nn(\alpha_{R_2/\<f\>})=\{\bv\in\sph_1:\bv\cdot\be_1\le0\}$
   and
   $\Nn(\alpha_{R_2/\<g\>})=\{\bv\in\sph_1:\bv\cdot\be_2\le0\}$,
   and so the intersection of these is the quarter-circle of
   $\sph_1$ in the third quadrant. However,
   \begin{displaymath}
      \{-\be_1,-\be_2\}=\Nn(\alpha_{R_2/\fa})
      \subsetneqq \Nn(\alpha_{R_2/\<f\>}) \cap
      \Nn(\alpha_{R_2/\<g\>}).
   \end{displaymath}
\end{example}

Given an ideal $\fa$ in $\rd$, we would like to compute a finite
generating set $F$ for $\fa$ with the property that if $H\in\Hd$
then $\fa$ contains an $H$-monic polynomial if and only if $F$ contains an
$H$-monic polynomial. For then
\begin{displaymath}
   \Nn(\arda)=\bigcap_{f\in F} \Nn(\ardf),
\end{displaymath}
where each term in the intersection has the description given
above.  

We compute such a set $F$ using the theory of Gr\"obner bases
\cite{AL} and universal Gr\"obner bases \cite{St}, which we now
briefly summarize. Let $\NN=\{0,1,2,\dots\}$. The monomials in
$\ZZ[u_1,\dots,u_d]$ correspond to elements $\bn\in\NN^d$ via
$u^{\bn}\leftrightarrow\bn$. A \textit{term order} $\prec$ on 
$\ZZ[u_1,\dots,u_d]$ is a total order on $\NN^d$ such that (1)
$\mathbf0\prec\bk$ for every $\mathbf0\ne\bk\in\NN^d$, and (2)
$\bk\prec\bm$ implies that $\bk+\bn\prec\bm+\bn$ for every $\bk$,
$\bm$, $\bn\in\NN^d$. Fix a term order $\prec$. For $0\ne f\in
\ZZ[u_1,\dots,u_d]$ let $\lt_{\prec}(f)$ denote the unique
\textit{leading term} $c_f(\bn)u^{\bn}$ of $f$, where $\bn$ is
maximal with respect to $\prec$ such that $c_f(\bn)\ne0$. Let
$\fb$ be an ideal of $\ZZ[u_1,\dots,u_d]$. Then a finite set
$G\subset\fa$ is a \textit{Gr\"obner basis for $\fb$ with respect
to $\prec$} if for every $f\in\fb$ there is a $g\in G$ for which
$\lt_{\prec}(g)$ divides $\lt_{\prec}(f)$ in
$\ZZ[u_1,\dots,u_d]$.

Although there are infinitely many term orders on
$\ZZ[u_1,\dots,u_d]$, it is shown in \cite[pp. 1--2]{St} that
every ideal has a \textit{universal Gr\"obner basis}, i.e., a
finite set that is a Gr\"obner basis with respect to
\textit{every} term order. Although the proof in \cite{St} is for
the case when the coefficients lie in a field, it is easy to adapt the
argument to integer coefficients. Also, \cite{St} provides
effective algorithms for computing a universal Gr\"obner basis.

\begin{proposition}\label{prop:groebner-basis}
   Let $\fa$ be an ideal in $\rd$. Then there is an effectively
   computable finite generating set $F\subset \fa$ such that
   \begin{equation}\label{eqn:grobner-intersection}
      \Nn(\arda)=\bigcap_{f\in F} \Nn(\ardf).
   \end{equation}
\end{proposition}

\begin{proof}
   The theory of Gr\"obner bases applies to polynomials rather
   than Laurent polynomials. In order to apply it here, we will
   subdivide $\zd$ into the $2^d$ orthants corresponding to the
   signs of the entries, and use a universal Gr\"obner basis in
   each of these orthants.

   Let $D=\{1,2,\dots,d\}$. For each subset $E\subset D$ define
   $E(j)=1$ if $j\in E$ and $E(j)=-1$ if $j\notin E$. Let
   $S_E=\ZZ[u_1^{E(1)},\dots,u_d^{E(d)}]$, and put $\fa_E=\fa\cap
   S_E$, an ideal in $S_E$. For each $E$ we will construct a
   finite set $F_E\subset\fa_E$, and show that
   $F=\bigcup_{E\subset D}F_E\subset\fa$ satisfies
   \eqref{eqn:grobner-intersection}.

   For notational simplicity we give the construction for $E=D$,
   and indicate the modifications needed for general $E$. Here
   $S_D=\ZZ[u_1,\dots,u_d]$. Choose a finite set
   $\{f_1,\dots,f_r\}$ of generators for $\fa_D$ over $S_D$. Let
   $R=S_D[t]=\ZZ[u_1,\dots,u_d,t]$ and define
   \begin{displaymath}
      \fb_D=\<f_1,\dots,f_r,t-u_1u_2\dots u_d\>\subset R
   \end{displaymath}
   (for general $E$ the last generator is $t-\prod_j
   u_j^{E(j)}$). Construct a universal Gr\"obner basis $G_D$ for
   $\fb_D$. Define $\phi\colon R\to\rd$ by $\phi(u_j)=u_j$ and
   $\phi(t)=u_1u_2\dots u_d$ (for general $E$ let $\phi(t)=\prod_j
   u_j^{E(j)}$). Clearly $\phi(\fb_D)=\fa_D$. Define
   $F_D=\phi(G_D)$.
   
   We show that if $\bv\ge0$ and $H=\Hv$, then there is an
   $H$-monic polynomial $f$ in $\fa$ if and only if some element
   in $F_D$ is also $H$-monic. (For general vectors $\bv$, use
   the above construction for $E=\{j:v_j\ge0\}$.) First, multiply
   $f$ by a monomial so that $f\in \fa_D$, and also that the
   $H$-exposed vertex of $f$ has the form $N\mathbf1$, where
   $\mathbf1=(1,1,\dots,1)$. Hence $f=u^{N\mathbf1}+h(u)$, where
   \begin{equation}\label{eqn:extremal-inequality}
      \bn\cdot\bv<(N\mathbf1)\cdot\bv \text{\quad for all
      $\bn\in\supp(h)$}.
   \end{equation}
   Therefore $t^N+h(u)\in\fb_D$.
   
   Define a term order on monomials in $R=S_D[t]$ by declaring
   that $(\bk,K)\prec(\bm,M)$ if and only if
   $(\bk+K\mathbf1)\cdot\bv\le(\bm+M\mathbf1)\cdot\bv$, and in
   case of equality that $K\ge M$, and in case of equality there,
   that $\bk$ is lexicographically less than $\bm$. Then
   $\lt_{\prec}(t^N+h(u))=t^N$ by
   \eqref{eqn:extremal-inequality}. Since $G_D$ is a Gr\"obner
   basis for $\prec$, there is a $g\in G_D$ such that $\lt_{\prec}(g)$
   divides $t^N$. Hence $\lt_{\prec}(g)=t^M$ for some $M\le N$.
   Let $\phi(g)=u^{\bm}+k(u)$ where $\phi(t^M)=u^{\bm}$.
   We claim that $\phi(g)$ is a polynomial in $F_D$ that is
   $H$-monic. This follows since
   every monomial in $g$ that could give rise to a term in
   $\phi(g)$ strictly larger than $\bm$ would have to
   already be greater than $t^M$ with respect to $\prec$.
\end{proof}

To compute $\N(\aM)$ for a general Noetherian $\rd$-module $M$,
we use the notion of Fitting ideal $\ff(M)$ to reduce to the
situation of Proposition ~\ref{prop:groebner-basis} (see
\cite[Chapter 20]{Eis} or \cite[\S XIII.10]{Lang} for
background). To define $\ff(M)$, suppose that $M$ is generated by
$m_1$, $m_2$, $\dots$, $m_r$. Let $K\subset\rd^r$ be the kernel
of the map $(f_1,\dots,f_r)\mapsto f_1m_1+\dots+f_rm_r$. Since
$K$ is Noetherian, it is finitely generated over $\rd$, say by
$(a_{11},\dots,a_{1r})$, $\dots$, $(a_{s1},\dots,a_{sr})$. Let
$A$ be the $s\times r$ matrix $[a_{ij}]$, so that $M\cong
\rd^r/(\rd^s A)$. Then $\ff(M)$ is defined to be the ideal in $\rd$
generated by all the $r\times r$ subdeterminants of $A$. This
ideal is independent of the presentation of $M$, and can be
effectively computed (see \cite{AL}). Also, if
$\ann(M)=\{f\in\rd:f\cdot M=0\}$, then
\begin{equation}\label{eqn:fitting}
   \bigl(\ann(M)\bigr)^r\subset\ff(M)\subset\ann(M)\subset
   \bigcap_{\fp\in\asc(M)}\fp.
\end{equation}

In \cite{EW} it is shown that both entropy and expansiveness of
$\aM$ can be computed from $\ff(M)$, although more subtle
dynamical information requires higher order Fitting ideals. The
next result shows that both pieces of $\N(\aM)$ can be found from~
$\ff(M)$.

\begin{proposition}\label{prop:fitting}
   Let $M$ be a Noetherian $\rd$-module and $\ff(M)$  denote
   its Fitting ideal. Then
   \begin{displaymath}
      \Nn(\aM)=\Nn(\ardfm), \quad \Nv(\aM)=\Nv(\ardfm),
   \end{displaymath}
   and so
   \begin{displaymath}
      \N(\aM)=\N(\ardfm).
   \end{displaymath}
\end{proposition}

\begin{proof}
   Suppose that $M$ is generated over $\rd$ by $r$ elements. Let
   $H\in\Hd$, and suppose that $M$ is $R_H$-Noetherian. By Lemma
   \ref{lem:noetherian-polynomial}, there is an $H$-monic
   polynomial $f\in\ann(M)$. Then \eqref{eqn:fitting} shows that
   $f^r\in\ff(M)$, and $f^r$ is clearly $H$-monic. Hence another
   application of Lemma \ref{lem:noetherian-polynomial} implies
   that $\rdfm$ is $H$-Noetherian. Thus
   $\Nn(\ardfm)\subset\Nn(\aM)$.
   Conversely, suppose that $\rdfm$ is $R_H$-Noetherian. Since
   $\ff(M)\subset\fp$ for every $\fp\in\asc(M)$, we see that
   $\rdp$ is $R_H$-Noetherian for every $\fp\in\asc(M)$. Hence
   \begin{displaymath}
      \Nn(\aM)=\bigcup_{\fp\in\asc(M)} \Nn(\ardp)\subset
      \Nn(\ardfm),
   \end{displaymath}
   proving that $\Nn(\aM)=\Nn(\ardfm)$.

   For the second equality, observe that since $\ff(M)\subset\fp$
   for every $\fp\in\asc(M)$,
   \begin{displaymath}
      \Nv(\aM)=\bigcup_{\fp\in\asc(M)}
      \Nv(\ardp)\subset\Nv(\ardfm). 
   \end{displaymath}
   To prove the reverse inclusion, take $\Hv\in\Nv(\ardfm)$. Then
   there is a $\bz\in\V(\ff(M))$ such that
   $\log|\bz|\in[0,\infty)\bv$. Suppose that $\bz\notin\V(\fp)$
   for every $\fp\in\asc(M)$. Then there are polynomials
   $g_{\fp}\in\fp$ with $g_{\fp}(\bz)\ne0$. Put
   $g=\prod_{\fp\in\asc(M)}g_{\fp}$. Then
   $g\in\bigcap_{\fp\in\asc(M)}\fp$, and so some power
   $g^k\in\ann(M)$ (to see this, either use a prime filtration of
   $M$, or observe that the radical of $\ann(M)$ is
   $\bigcap_{\fp\in\asc(M)}\fp$). By \eqref{eqn:fitting} we see
   that $(g^k)^r\in\ff(M)$. But this contradicts
   $g^{kr}(\bz)=\prod_{\fp} g_{\fp}^{kr}(\bz)\ne0$. Hence there
   is a $\fp\in\asc(M)$ such that $\Hv\in\Nv(\ardp)$.
\end{proof}

\begin{remark}\label{rem:bieri}
   The set $\Nn(\al_M)$ has been investigated in another context
   by Bieri and Groves \cite{BG} using a valuation-theoretic
   approach to the extension of characters on fields.  More
   specifically, they show in \cite[Theorem~8.1]{BG} how $\Nn(\aM)$
   can be calculated using such characters. In \cite{Miles} these
   ideas are developed to give a purely valuation-theoretic
   description of all of $\N(\aM)$.
\end{remark}

Next, we characterize when $\Nn(\al_M)$ or $\Nv(\al_M)$ is empty.
Note that the answers depend only on the topological nature of
$X_M$. If $\fp$ is a prime ideal in $\rd$, we let
$\charact(\rdp)$ denote the characteristic of the integral domain
$\rdp$. 

\begin{proposition}\label{prop:empty}
   Let $M$ be a Noetherian $\rd$-module.
   \begin{enumerate}
     \item[(1)] The following conditions are equivalent.
      \begin{enumerate}
        \item[(a)] $\Nn(\al_M)=\emptyset$.
        \item[(b)] $M$ is finitely generated as an abelian group.
        \item[(c)] $X_M$ is the direct product of a
         finite-dimensional torus and a finite abelian group.
      \end{enumerate}
     \item[(2)] The following conditions are equivalent.
      \begin{enumerate}
        \item[(a)] $\Nv(\al_M)=\emptyset$.
        \item[(b)] $M$ is a torsion abelian group.
        \item[(c)] $X_M$ is totally disconnected.
      \end{enumerate}
   \end{enumerate}
\end{proposition}

\begin{proof}
   (1b)\implies(1a): Suppose that $M$ is finitely generated as an
   abelian group. Then $M$ is trivially $R_H$-Noetherian for
   every $H\in\Hd$, so that $\Nn(\al_M)=\emptyset$ by Remark~
   \ref{rem:non-Noetherian-M}.

   (1a)\implies(1b): Our proof uses an algebraic analogue of the
   proof of Theorem 3.6 in \cite{BL}.

   Since $\Nn(\al_M)=\emptyset$, by Lemma
   \ref{lem:noetherian-polynomial} for every $\bv\in\sphd$ there
   is a polynomial of the form $1-f(u)=1-\sum_{\bn\in
   F}c_f(\bn)u^{\bn}\in R_{\Hv}$ that annihilates $M$ and such
   that $\bn\cdot\bv<0$ for every $\bn\in F$. This shows that
   there are $\veps_{\bv}>0$, $r_{\bv}>0$, and a neighborhood
   $\SU(\bv)$ in $\sphd$ such that if $\bp\in\RR^d$ with
   $\|\bp\|>r_{\bv}$ and $\bp/\|\bp\|\in\SU(\bv)$, then
   $F+\bp\subset B(\|\bp\|-\veps_{\bv})$.

   The collection $\{\SU(\bv):\bv\in\sphd\}$ is an open cover of
   $\sphd$, so by compactness there is a finite subcover
   $\{\SU(\bv_1),\dots,\SU(\bv_n)\}$. Put
   $\veps=\min\{\veps_{\bv_j}:1\le j\le n\}$ and
   $r=\max\{r_{\bv_j}:1\le j\le n\}$. Let $M$ be generated over
   $\rd$ by $m_1$, $\dots$, $m_k$. For $s>0$ let $M_s$ denote the
   abelian group generated by $\{u^{\bn}\cdot m_i:1\le i\le k,
   \bn\in B(s)\cap\zd\}$. We claim that $M_r=M$, so that $M$ is
   finitely generated as an abelian group. We prove this by
   showing successively that
   $M_r=M_{r+\veps}=M_{r+2\veps}=\cdots$, so that
   $M=\bigcup_{q=0}^\infty M_{r+q\veps}=M_r$.

   Let $\bp\in\bigl( B(r+\veps)\smallsetminus
   B(r)\bigr)\cap\zd$. Then $\bp/\|\bp\|\in\SU(\bv_j)$ for some
   $1\le j \le n$. By our construction, there is a polynomial
   $1-f(u)=1-\sum_{\bn\in F}c_f(\bn)u^{\bn}$ that annihilates
   $M$, and such that $F+\bp\subset
   B(\|\bp\|-\veps_{\bv_j})\subset B(r)$. Hence for $1\le i\le k$
   we see that
   \begin{displaymath}
      u^{\bp}\cdot m_i=u^{\bp}f(u)\cdot m_i=\sum_{\bn\in F}
         c_f(\bn)u^{\bn+\bp}\cdot m_i\in M_r.
   \end{displaymath}
   This shows that $M_{r+\veps}\subset M_r$, and the reverse
   inclusion is trivial.  The same argument applied to
   $M_{r+\veps}$ shows that $M_{r+\veps}=M_{r+2\veps}$, and so
   on, completing the proof.

   (1b)\ifff(1c): This equivalence is standard from duality.

   (2a)\implies(2b): By assumption, $\Nv(\ardp)=\emptyset$ for
   every $\fp\in\asc(M)$. Hence $\charact(\rd/\fp)>0$ for every
   $\fp\in\asc(M)$ by Hilbert's Nullstellensatz. There is a prime
   filtration $0=M_0\subset M_1\subset \dots \subset M_r=M$ with
   $M_j/M_{j-1}\cong \rd/\fq_j$, where $\fq_j$ is a prime ideal
   containing some $\fp_j\in\asc(M)$. Then the integer $\prod_{j=1}^r
   \charact(\rd/\fp_j)$ annihilates $M$, so that $M$ is a torsion
   abelian group.

   (2b)\implies(2a): Suppose that $M$ is a torsion abelian group.
   Let $\fp\in\asc(M)$. Then there is an $m\in M$ with $\rd\cdot
   m\cong\rd/\fp$. Since $\rd\cdot m$ must also be a torsion
   abelian group, it follows that $\charact(\rd/\fp)>0$. In
   particular, $\fp$ contains a nonzero constant, so that
   $\Nv(\ardp)=\emptyset$. Hence $\Nv(\al_M)=\emptyset$.

   (2b)\ifff(2c): This equivalence is standard from duality.
\end{proof}

\section{Expansive rank, entropy rank, and Krull dimension}
   \label{sec:ranks}

Although a $\zd$-action nominally involves $d$ commuting
transformations, sometimes its true ``rank'' is less.  In this
section we describe two ways to measure this rank. In what
follows $\entropy$ denotes topological entropy.

\begin{definition}\label{def:ranks}
   Let $\beta$ be a topological $\zd$-action. Define the
   \textit{expansive rank} of $\beta$ to be
   \begin{displaymath}
      \exprk(\beta)=\min\{ k:\E_k(\beta)\ne\emptyset \},
   \end{displaymath}
   and the \textit{entropy rank} of $\beta$ to be
   \begin{displaymath}
      \entrk(\beta)=\max\{ k:\text{there is a rational $k$-plane
      $V$ with $\entropy(\beta,V\cap\zd)>0$} \}.
   \end{displaymath}
   By convention, if $\beta$ is not expansive we put
   $\exprk(\beta)=d+1$, and if the set defining $\entrk(\beta)$
   is empty we put $\entrk(\beta)=0$. 
\end{definition}

These ranks attempt to measure, from the viewpoints of
expansiveness and of entropy, the maximum number of ``independent
transformations'' in the action such that the remaining
transformations are determined in some sense. As a concrete
instance of what we have in mind, consider Ledrappier's Example
\ref{exam:ledrappier}. If $L$ is an expansive line, and $L^1$
denotes the thickening of $L$ by $1$, then every element of the
action can be written as a function of $\al^{\bn}$ for finitely
many $\bn\in L^1\cap\ZZ^2$. In this sense $\al$ has expansive
rank $1$.

\begin{proposition}\label{prop:entrk<exprk}
   Let $\beta$ be a topological $\zd$-action. Then
   $\entrk(\beta)\le\exprk(\beta)$. 
\end{proposition}

\begin{proof}
   Let $\exprk(\beta)=k$. It follows from \cite[Theorem 6.3]{BL}
   that if $V$ is a $k$-plane, then $\entropy(\beta,V\cap\zd)<\infty$. A
   standard argument now shows that if $W$ is a rational subspace of
   dimension $\ge k+1$ then $\entropy(\beta,W\cap\zd)=0$. Hence
   $\entrk(\beta)\le k$.
\end{proof}

We remark that an alternative proof of this proposition can be
obtained by modifying arguments of Shereshevsky \cite{Sher}.
We also observe that in general the inequality here can be strict,
for example a zero entropy subshift has entropy rank zero
but expansive rank one. 

The building blocks for algebraic $\zd$-actions are based on
modules of the form $\rdp$ for prime ideals $\fp$. For these we
investigate expansive and entropy ranks. The appropriate
algebraic version of rank is Krull dimension for rings.  Recall
that the \textit{Krull dimension} $\kdim \mathscr{R}$ of a ring
$\mathscr{R}$ is the supremum of the lengths $r$ of all chains
$\fp_0\subsetneqq\fp_1\subsetneqq\dots\subsetneqq\fp_r$ of prime
ideals in $\mathscr{R}$ (see \cite[Chapter 8]{Eis} for the
necessary background). For example, $\kdim R_d=d+1$. In
\cite[Theorem 7.5]{BL} it is shown that if $\fp$ is a prime ideal
in $R_d$ generated by $g$ elements, then
\begin{displaymath}
   \exprk(\ardp)\ge \kdim(\rdp)-1\ge d-g,
\end{displaymath}
and if $\charact(\rdp)>0$ then
\begin{displaymath}
   \exprk(\ardp)\ge d-g+1.
\end{displaymath}
Since expansive and entropy ranks in some sense measure
dimension, one might heuristically expect $\exprk(\ardp)$,
$\entrk(\ardp)$ and $\kdim(\rdp)$ to coincide. The best that can
be said in this direction is as follows.

\begin{proposition}\label{prop:ranks}
   Let $\ardp$ be an expansive action with zero entropy. Then
   \begin{displaymath}
      \entrk(\ardp)=\kdim(\rdp)\le\exprk(\ardp).
   \end{displaymath}
   If $\charact(\rdp)>0$ then all of these quantities are equal.
\end{proposition}

\begin{proof}
   Let $k=\exprk(\ardp)$. Since the expansive set is open, there
   is a rank $k$ lattice in $\zd$ for which the restriction of
   $\ardp$ is expansive. By a simple change of variables, we may
   assume that this lattice is generated by the first $k$ unit
   vectors. That is, the subaction dual to multiplication by
   $u_1$, $\dots$, $u_k$ is expansive. We thus consider $\rdp$ as
   an $R_k$-module. We claim that the images of the variables
   $u_1$, $\dots$, $u_k$ in $\rdp$ must satisfy a polynomial
   relation, and also that each $u_j$ for $j>k$ must be algebraic
   over $u_1$, $\dots$, $u_k$. For the first claim, observe that
   if the images of $u_1$, $\dots$, $u_k$ in $\rdp$ do not
   satisfy any polynomial relation, then the natural map
   $R_k\to\rdp$ is injective, so that $\{0\}$ is a prime ideal
   associated to the $R_k$-module $\rdp$. But by
   Theorem~\ref{thm:expansive} this contradicts the assumption
   that the restriction of $\ardp$ to the first $k$ variables is
   expansive. If some $u_j$ were not algebraic over $u_1$,
   $\dots$, $u_k$, then $\rdp$ would not be a Noetherian
   $R_k$-module, again contradicting the expansiveness
   assumption. It follows that $\kdim(\rdp)\le \kdim(R_k)-1=k$.
   
   Now let $\kdim(\rdp)=k$, and assume first that
   $\charact(\rdp)=0$. Then any $k+1$ monomials must satisfy two
   coprime irreducible polynomial relations, so that the entropy
   of the corresponding $\ZZ^{k+1}$-action is zero. Hence
   $\entrk(\ardp)\le\kdim(\rdp)$. Now assume without loss of
   generality that the images of the variables $u_{k+1}$,
   $\dots$, $u_d$ in $\rdp$ are algebraic over the images of
   $u_1$, $\dots$, $u_k$, and that $u_1$, $\dots$, $u_k$ satisfy
   only one polynomial relation. If this relation is not a
   generalized cyclotomic polynomial, then the corresponding
   $\ZZ^k$-subaction has positive entropy \cite[Example
   5.4]{LSW}, so that $\entrk(\ardp)\ge\kdim(\rdp)$. So we may
   assume that this single polynomial relation is an irreducible
   generalized cyclotomic polynomial. After a suitable change of
   variables, this implies that $u_k^r=1$ in $\rdp$ for some
   $r\ge1$. Since the original system is expansive, we may find
   among the variables $u_{k+1}$, $\dots$, $u_d$ a variable $u_j$
   which is not a root of unity. Now the same argument may be
   applied to the set of variables $u_1$, $\dots$, $u_{k-1}$,
   $u_j$. Continuing, we either arrive at a $\ZZ^k$-subaction
   with positive entropy or a contradiction. We deduce that
   $\entrk(\ardp)\ge\kdim(\rdp)$, and so
   $\entrk(\ardp)=\kdim(\rdp)$.

   Finally, suppose $\charact(\rdp)=p>0$, and let
   $k=\kdim(\rdp)$. Then $\rdp$ is a ring extension of
   $\mathbb{F}_p[u_1^{\pm1},\dots, u_k^{\pm1}]$. By Noether
   normalization \cite[Section~ 8]{Sch} we may choose variables so
   that $u_1$, $\dots$, $u_k$ do not satisfy a polynomial
   relation, and the variables $u_{k+1}^{\pm1}$, $\dots$,
   $u_d^{\pm1}$ are integral over $u_1$, $\dots$, $u_k$. This
   shows that there is an expansive, positive-entropy
   $\ZZ^k$-subaction. Hence $\entrk(\ardp)=\exprk(\ardp)=\kdim(\rdp)$.
\end{proof}

The next example shows that the inequality in Proposition
\ref{prop:ranks} can be strict when $\charact(\rdp)=0$.

\begin{example}\label{exam:strict-inequality}
   (\textit{Krull dimension strictly less than expansive rank}.)
   Let $\phi(z)=z^2-z-i=(z-\lambda_1)(z-\lambda_2)$ and
   $\psi(z)=z^2+z-2i=(z-\mu_1)(z-\mu_2)$, where $i=\sqrt{-1}$.
   Set $F=\{0,\lambda_1,\lambda_2,\mu_1,\mu_2\}$. Define
   $\bs\colon\CC\to\CC^3$ by
   $\bs(z)=\bigl(z,\phi(z),\psi(z)\bigr)$. Put
   \begin{displaymath}
      W=\{\,\bs(z):z\in\CC\smallsetminus F\,\}\subset(\CC^\times)^3.
   \end{displaymath}
   Let $\fp\subset R_3$ be the ideal of Laurent polynomials
   in $R_3$ that vanish on all of $W$.

   We will show that $\fp$ is a prime ideal, that $\arthp$ is
   expansive and mixing, and that
   $2=\kdim(\rthp)<\exprk(\arthp)=3$.

   To begin, let $g(u,v,w)=(u^2-u-v)^2+1$ and
   $h(u,v,w)=w-2v+u^2-3u$. We will show that $\fp=\<
   g,h\>$. Both $g$ and $h$ vanish on $W$, and so $\<
   g,h\>\subset\fp$. To establish the reverse inclusion,
   suppose that $f\in\fp$. We can reduce $f$ in $R_3$ modulo
   $\<h\>$ to have the form $w^{-n}k(u,v)$ for some $n\in\ZZ$ and
   $k\in R_2$. We can then reduce $k(u,v)$ modulo $\<g\>$ in
   $R_2$ to have the form $v^{-m}[p(u)v+q(u)]$ for some $m\in\ZZ$
   and $p(u),q(u)\in R_1$. Hence $f$ is congruent in $R_3$ to
   $v^{-m}w^{-n}[p(u)v+q(u)]$ modulo $\<g,h\>$. Since $f$ is in
   $\fp$, it vanishes on $W$, so that
   \begin{displaymath}
      f(\bs(z))=\phi(z)^{-m}\psi(z)^{-n}\bigl
      [ p(z)\phi(z)+q(z)\bigr]=0
      \text{\quad for $z\in\CC\smallsetminus F$.}
   \end{displaymath}
   Hence as Laurent polynomials with complex coefficients,
   $p(z)\phi(z)=-q(z)$. If $p\ne0$ in $R_1$, then the nonzero
   coefficient of the smallest power of $z$ on the left-hand side
   would be pure imaginary, while on the right-hand side it would
   be real. This contradiction shows that $p=0$, and thus $q=0$,
   hence $f\in\<g,h\>$. This proves that $\fp=\<g,h\>$.

   To show that $\fp$ is prime, suppose that $f\cdot k\in\fp$ for
   some $f,k\in R_3$. Then $f(\bs(z))\cdot k(\bs(z))$ is  a
   rational function that vanishes on $\CC\smallsetminus
   F$. Hence at least one of $f(\bs(z))$ or $k(\bs(z))$ must
   vanish on $\CC\smallsetminus F$, so that $f\in\fp$ or
   $k\in\fp$, establishing primality of $\fp$. Alternatively, one
   can give an algebraic proof that $\fp$ is prime along the
   lines of Example \ref{exam:eins-ward}.

   To compute $\kdim(\rthp)$, first note that $h$ gives $w$ in
   terms of $u$ and $v$ in $\rthp$. Then $g$ shows that $v$ is
   algebraic over $u$ in $\rthp$. It now follows from standard
   facts about Krull dimension (see \cite[Chapter 13]{Eis}) that
   $\kdim(\rthp)=\kdim(R_1)=2$.
   
   We next examine the expansive character of $\arthp$. Define
   $\overline{W}=\{(\overline{a},\overline{b},\overline{c}):(a,b,c)\in
   W\}$. We claim that $\V(\fp)=W\cup\overline{W}$. Since all
   polynomials in $\fp$ have real coefficients, it is clear that
   $W\cup\overline{W}\subset\V(\fp)$. Conversely, suppose that
   $(a,b,c)\in\V(\fp)$. Since $g(a,b,c)=0$, it follows that
   $a^2-a-b=\pm i$. If $a^2-a-b=i$, then using $h(a,b,c)=0$ we
   see that $c=a^2+a-2i$. Thus $(a,b,c)=\bs(a)\in W$. If
   $a^2-a-b=-i$, then taking complex conjugates and applying the
   previous argument shows that
   $(\overline{a},\overline{b},\overline{c})
   =\bs(\overline{a})\in W$, so that
   $(a,b,c)\in\overline{W}$. Hence $\V(\fp)=W\cup\overline{W}$.
   Note that $\log|W|=\log|\overline{W}|$, so that
   $\log|\V(\fp)|=\log|W|$.
   
   By Theorem \ref{thm:expansive}, $\arthp$ is expansive if and
   only if $\mathbf{0}\notin\log|W|$. Hence it suffices to show
   that the curve
   $\gamma(\theta)=(\log|\phi(e^{i\theta})|,\log|\psi(e^{i\theta})|)$
   does not pass through the origin. But this is clear from the graph
   of $\gamma$ shown in Figure~ \ref{fig:amoeba}(a). Hence
   $\arthp$ is expansive.

   \begin{figure}[htbp]
      \begin{center}
          \scalebox{1.0}{\includegraphics{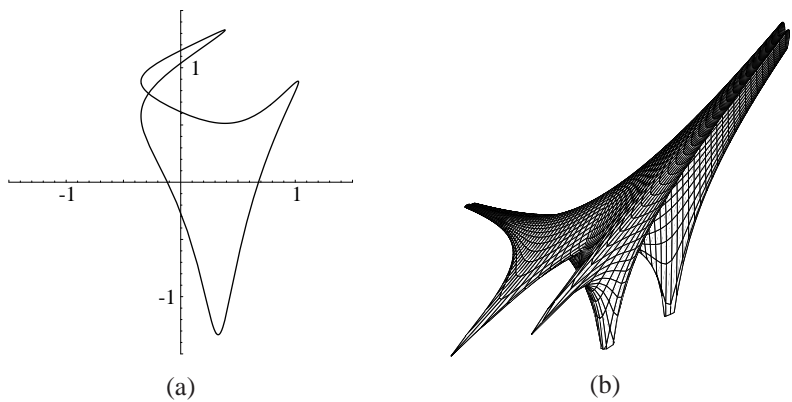}}
         \caption{Cross-section and logarithmic image}
         \label{fig:amoeba}
      \end{center}
   \end{figure}
   
   Next, we show that $\Nv(\arthp)=\sph_2$, and hence that
   $\N(\arthp)=\sph_2$. Consider that
   cross-section of $\log|W|$ corresponding to setting the first
   coordinate to $\log r$.  Since then $z=re^{i\theta}$, this
   cross-section is given by the curve
   $\gamma_r(\theta)=\bigl(\log|\phi(re^{i\theta})|,
   \log|\psi(re^{i\theta})|\bigr)$.  For $r$ close to $1$ each
   $\gamma_r$ is a curve that surrounds the origin, and hence the
   part of $\log|W|$ with first coordinate close to $0$ is a
   surface that surrounds the origin in $\RR^3$. As $r$ varies
   from $0$ to $\infty$, this surface generates six ``tendrils''
   stretching off to infinity (see Figure~\ref{fig:amoeba}(b),
   where the term ``amoeba'' is apt). Two tendrils correspond to
   $r=|\lambda_1|\approx 1.44$ and $r=|\lambda_2|\approx 0.69$,
   and they are asymptotic to the half-lines $L_1$ and $L_2$ in
   the $-\be_2$ direction given by
   $(\log|\lambda_j|,0,\log|\psi(\lambda_j)|)-[0,\infty)\be_2$
   for $j=1,2$. Another two tendrils correspond to the half-lines
   $L_3$ and $L_4$ in the $-\be_3$ direction given by
   $(\log|\mu_k|,\log|\phi(\mu_k)|,0)-[0,\infty)\be_3$ for
   $k=1,2$. As $r\to0$, $\log r\to -\infty$, and there is a fifth
   tendril asymptotic to the half-line $L_5$ in the $-\be_1$
   direction given by
   $(0,\log|\phi(0)|,\log|\psi(0)|)-[0,\infty)\be_1$. Finally, as
   $r\to\infty$ there is a sixth tendril asymptotic to
   $L_6=[0,\infty)\bv$, where $\bv=(1/3,2/3,2/3)$.

   Note that when the half-lines $L_1$ through $L_5$ are extended
   to lines, none passes through the origin. It follows that the
   radial projection of $\log|W|$ to $\sph_2$ must cover all of
   $\sph_2$, with the possible exception of $\bv=(1/3,2/3,2/3)$
   corresponding to $L_6$. At this point we can conclude that
   $\N(\arthp)=\sph_2$ since it is closed. However, we continue
   with the complete description of $\Nv(\arthp)$.

   To handle the remaining point $\bv$, we show that there is a
   $z\in\CC\smallsetminus F$ such that
   $\log|\phi(z)|=\log|\psi(z)|=2\log|z|$. The equation
   $\log|\phi(z)|=2\log|z|=\log|z^2|$ implies that
   $|\phi(z)/z^2|=1$, so that $\phi(z)/z^2=e^{i\theta}$ for some
   $\theta$. Hence $(1-e^{i\theta})z^2-z-i=0$, with roots
   \begin{equation}\label{eqn:first-root}
      z=\frac{ 1\pm\sqrt{ 1+4i(1-e^{i\theta}) } }
            { 2(1-e^{i\theta}) }.
   \end{equation}
   A similar calculation starting with $\log|\psi(z)|=2\log|z|$
   shows that
   \begin{equation}\label{eqn:second-root}
      z=\frac{ -1\pm\sqrt{ 1+8i(1-e^{i\xi}) } }
            { 2(1-e^{i\xi}) }
   \end{equation}
   for some $\xi$. Plotting the roots using the positive sign for
   \eqref{eqn:first-root} and \eqref{eqn:second-root} gives the
   graphs shown in Figure \ref{fig:graph}. They cross in two
   points, one of which is $z_1\approx 0.53+3.36i$. Then
   $z_1\in\CC\smallsetminus F$, $\log|z_1|>0$,
   $\log|\phi(z_1)|>0$, $\log|\psi(z_1)|>0$, and so
   $\log|\bs(z_1)|=t\bv$ for some $t>0$. Hence
   $\bv\in\Nv(\arthp)$, and so $\Nv(\arthp)=\sph_2$. Thus
   $\exprk(\arthp)=3$, showing that here Krull dimension is
   strictly less than expansive rank.

   \begin{figure}[htbp]
      \begin{center}
          \scalebox{1.0}{\includegraphics{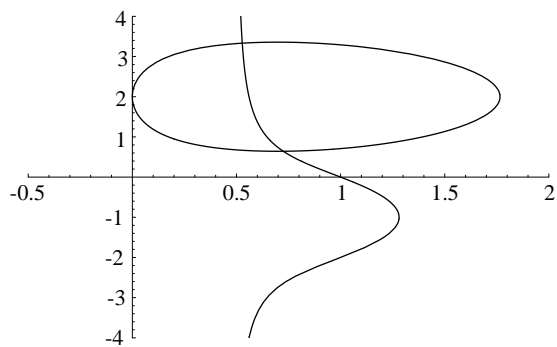}}
         \caption{Graphs of roots}
         \label{fig:graph}
      \end{center}
   \end{figure}

   We conclude this example with two further properties of
   $\arthp$: it is mixing, and every 2-dimensional plane has
   positive entropy.

   By Theorem 6.5(2) of \cite{Sch}, $\arthp$ is mixing if and
   only if $u^{\bn}-1\notin\fp$ for all
   $\bn\in\ZZ\smallsetminus\{\mathbf{0}\}$. Suppose the contrary.
   Then the plane $\{\bx\in\RR^3:\bx\cdot\bn=0\}$ would contain
   $\log|W|$, which is evidently false. Thus $\arthp$ is mixing.
   
   Let $Q$ be a rational 2-dimensional subspace of $\RR^3$.
   Denote the lattice $Q\cap\ZZ^3$ by $\Gamma$, and let
   $\bm,\bn\in\Gamma$ form an integral basis for $\Gamma$.
   Suppose that $\entropy(\arthp,\Gamma)=0$. Let
   $\rg=\ZZ[u^{\bk}:\bk\in\Gamma]$ and $\pg=\fp\cap\rg$. As
   $\Gamma$-actions, $\alpha_{\rg/\pg}$ is a factor of $\arthp$,
   so that $\entropy(\alpha_{\rg/\pg})=0$. By the preceding paragraph,
   $\pg$ contains no generalized cyclotomic polynomials, so by
   Theorem 4.2 of \cite{LSW}, it follows that $\pg$ cannot be
   principal. Hence $\kdim(\rg/\pg)=1$ and $\rg/\pg$ is certainly
   $\ZZ$-torsion-free, so that $\kdim[(\rg/\pg)\otimes\QQ]=0$.
   Thus $\V(\pg)$ is finite, say $\V(\pg)=\{(\xi_j,\eta_j):1\le
   j\le r\}$. Hence for every $\bz\in\V(\fp)$ there is a $1\le j
   \le r$ for which $\bz^{\bm}=\xi_j$ and $\bz^{\bn}=\eta_j$, or
   $\bm\cdot\log|\bz|=\log|\xi_j|$ and
   $\bn\cdot\log|\bz|=\log|\eta_j|$. It follows that if $L_j$
   denotes the line of intersection between the two planes
   $\{\bx\in\RR^3:\bm\cdot\bx=\log|\xi_j|\}$ and $\{\bx\in\RR^3:
   \bn\cdot\bx=\log|\eta_j|\}$, then
   $\log|\V(\fp)|\subset\bigcup_{j=1}^r L_j$. But the radial
   projection of finitely many lines to $\sph_2$ cannot cover
   $\sph_2$, contradicting $\Nv(\arthp)=\sph_2$. This shows that
   $\arthp$ has strictly positive entropy with respect to every
   rational 2-dimensional plane.
\end{example}

Example \ref{exam:strict-inequality} shows that 
when $\kdim(\rdp)=d-1$, then $\Nv(\ardp)$ can
have nonempty interior. We now turn
to $\Nn(\ardp)$. We have already observed two cases
when $\Nn(\ardp)$ has nonempty interior: $\fp=0$, with
$\kdim(\rdp)=d+1$; and $\fp=\<f\>$ with $|c_f(\bn)|>1$ for some
vertex $\bn$ of $\SN(f)$ (see Example
\ref{exam:three-variables-principal}), where 
$\kdim(\rdp)=d$. Our next result shows that when $\kdim(\rdp)$
drops below ~$d$, then $\Nn(\ardp)$ cannot have interior.

\begin{proposition}\label{prop:Nn-empty-interior}
   Suppose that $\fp$ is a prime ideal in $\rd$ with
   $\kdim(\rdp)\le d-1$. Then $\Nn(\ardp)$ is closed with empty
   interior. 
\end{proposition}

\begin{proof}
   Remark \ref{rem:Nn-is-closed} shows that $\Nn(\ardp)$ is
   closed.

   First suppose that $\charact(\rdp)=p>0$. Let
   $R_{d,p}=(\ZZ/p\ZZ)[u_1^{\pm1},\dots,u_d^{\pm1}]$, and $\fq$
   be the image of $\fp$ under the map $R_d \to R_{d,p}$ that
   reduces coefficients mod~ $p$. Then $R_d/\fp\cong
   R_{d,p}/\fq$. Hence
   $\kdim(R_{d,p}/\fq)=\kdim(R_d/\fp)\le
   d-1<d=\kdim(R_{d,p})$, so that $\fq\ne0$. Choose $0\ne
   f\in\fq$. Then $\Nn(\ardp)$ is contained in
   $\Nn(\alpha_{R_{d,p}/\<f\>})$, which by \cite[Theorem~
   7.2]{BL} equals the $(d-2)$-skeleton of the spherical dual of
   the mod~ $p$ Newton polyhedron of $f$. Hence $\Nn(\ardp)$ has
   empty interior.

   Finally, suppose that $\charact(\rdp)=0$. Recall from
   Definition \ref{def:H-monic} the meaning of support and of
   $H$-exposed vertex.
   
   Let $\SU$ be an arbitrary nonempty open set in $\Hd$. Choose
   $a\ge1$ minimal so that there is an $H\in\SU$ and $f\in\fp$
   having an $H$-exposed vertex $\bn$ with $c_f(\bn)=a$. Note
   that $\fp\ne0$ since $\kdim(\rdp)\le d-1$, so that $a$
   exists. We claim that $a=1$. If so, then Lemma
   \ref{lem:noetherian-polynomial} shows that $\SU$ contains an
   $H$ for which $\rdp$ is $R_H$-Noetherian, and so $\Nn(\ardp)$
   is nowhere dense.

   Suppose that $a>1$. Choose $H\in\SU$ and $f\in\fp$ such that
   $f$ has an $H$-exposed vertex $\bn$ with $c_f(\bn)=a$. We may
   clearly assume that $f$ is irreducible since any factor of $f$
   would also have an $H$-exposed vertex. Since $H$-exposure is
   an open condition and rational directions are dense, we
   may also assume that $H$ is rational, i.e., that there is an
   $\bm\in\zd\smallsetminus\{\mathbf{0}\}$ with $H=\{\bx\in
   \RR^d:\bx\cdot\bm\le0\}$. Adjusting $f$ by a monomial if
   necessary, we may assume that $\supp(f)\subset H$ and that
   $\supp(f)\cap \partial H=\{\mathbf0\}$, so that
   $c_f(\mathbf0)=a$.

   Next, we claim that there is a $0\ne g \in\fp$ with
   $\supp(g)\subset\partial H$. Equivalently, if $R_{\bm}$
   denotes $\ZZ[u^{\bk}:\bk\cdot\bm=0]$, then
   $R_{\bm}\cap\fp\ne\{0\}$. For suppose that
   $R_{\bm}\cap\fp=\{0\}$.
   Since $f$ is irreducible and $\rd/\<f\>$ has Krull dimension $d$,
   it follows that $0\subsetneqq
   \<f\>\subsetneqq\fp$ is a chain of prime ideals in $\rd$. Let
   $S$ denote the multiplicatively closed subset
   $R_{\bm}\smallsetminus\{0\}$ of~ $\rd$. Since
   $S\cap\fp=\emptyset$, it follows that $0\subsetneqq
   S^{-1}(\<f\>)\subsetneqq S^{-1}\fp$ is a chain of prime ideals
   in $S^{-1}\rd$. But $S^{-1}\rd\cong S^{-1}R_{\bm}[v,v^{-1}]$
   for a suitable monomial $v$,
   and $S^{-1}R_{\bm}$ is a field, so that $\kdim\bigl(
   S^{-1}R_{\bm}[v,v^{-1}]\bigr)=1$, contradicting the existence of a
   chain of primes of length two. Hence $R_{\bm}\cap\fp\ne\{0\}$.

   Let $0\ne g\in R_{\bm}\cap\fp$. Since $\charact(\rdp)=0$, then
   $\fp\cap\ZZ=\{0\}$. If all the coefficients of $g$ are
   divisible by $a$, then primality of $\fp$ shows that $g/a$ is
   also in $\fp$. Hence we can find $g\in R_{\bm}\cap\fp$ not all
   of whose coefficients are divisible by $a$. Let
   \begin{displaymath}
      h=g\;-\sum_{\bn\in\supp(g)} \Bigl\lfloor
      \frac{c_g(\bn)}{a}\Bigr\rfloor
      \,u^{\bn}f.
   \end{displaymath}
   Since $\supp(g)\subset\partial H$, $\supp(f)\cap\partial
   H=\{\mathbf0\}$, and $c_f(\mathbf0)=a$, it follows that
   $h\in\fp$, $\supp(h)\subset H$, $\supp(h)\cap\partial
   H\ne\emptyset$, and $0<c_h(\bn)<a$ for every
   $\bn\in\supp(h)\cap\partial H$. Hence there is a small
   perturbation $H'\in\SU$ of $H$ for which $h$ has an
   $H'$-exposed vertex $\bn$ with $0<c_h(\bn)<a$. This
   contradicts minimality of $a$, proving that $a=1$, and
   completing the proof.
\end{proof}

\section{Lower-dimensional subspaces}\label{sec:lower-dimensional}

Thus far we have concentrated on expansive behavior along
half-spaces and their $(d-1)$-dimensional boundaries. Once this
is found, then expansive behavior along lower-dimensional
subspaces is completely determined by \cite[Theorem 3.6]{BL}: a
$k$-plane is nonexpansive for a $\zd$-action if and only if it is
contained in a nonexpansive $(d-1)$-plane (or, equivalently, a
nonexpansive half-space). Additionally, half-spaces $H$ give rise
to subrings $R_H$ of $\rd$, and this algebraic structure makes
certain arguments work smoothly.

In this section we sketch how to modify our definitions and
proofs to work for lower-dimensional subspaces. In particular, we
obtain a direct description of the set $\N_k(\aM)$ of
nonexpansive $k$-planes as the union of two pieces, one coming
from a Noetherian condition on $M$ along $k$-planes, and the
other from a variety condition involving the orthogonal
complements of $k$-planes. This description is a
lower-dimensional version of \eqref{eqn:N-as-union}.

We begin by defining a notion of Noetherian along a general
subspace of $\RR^d$.

\begin{definition}\label{def:Noetherian-along-V}
   Let $M$ be a Noetherian $\rd$-module and $V$ be a $k$-plane in
   $\RR^d$. Then $M$ is \textit{Noetherian along~ $V$} if $M$ is
   $R_H$-Noetherian for every half-space $H\in\Hd$ containing
   $V$. The set of $k$-planes $V$ along which $M$ is not
   Noetherian is denoted by $\Nn_k(\al_M)$.
\end{definition}

\begin{remark}\label{rem:Nnk-closed}
   By Remark \ref{rem:Nn-is-closed}, $\{H\in\Hd:\text{$M$ is
   $R_H$-Noetherian}\}$ is open. It follows that $\{V\in\G_k:
   \text{$M$ is Noetherian along~ $V$}\}$ is open, so that
   $\Nn_k(\al_M)$ is closed in $\G_k$.
\end{remark}

In the following, recall that $V^t$ denotes the thickening of a
subspace $V$ by an amount $t$.

\begin{lemma}\label{lem:Noetherian-along-V}
   Let $M$ be a Noetherian $\rd$-module and $V$ be a subspace of
   $\RR^d$. Then $M$ is Noetherian along~ $V$ if and only if there
   are $m_1$, $\dots$, $m_r\in M$ and $t>0$ such that $M$ is
   generated as a group by $\{u^{\bn}m_j:1\le j\le r \text{ and }
   \bn\in V^t\cap\zd\}$. In particular, if $V$ is rational, then
   $M$ is Noetherian along~ $V$ if and only if $M$ is Noetherian
   as a module over the ring $\ZZ[u^{\bn}:\bn\in V\cap\zd\}$.
\end{lemma}

\begin{proof}
   The case $V=0$ is the implication (1a)\implies(1b) in
   Proposition \ref{prop:empty}. The case of general $V$ uses an
   entirely analogous adaptation of the proof of \cite[Theorem
   3.6]{BL}.

   If $V$ is rational, let $R_V=\ZZ[u^{\bn}:\bn\in
   V\cap\zd]$. Then $R_V$ is a Noetherian ring. If $M$ is
   generated as a group by $\{u^{\bn}m_j:1\le j\le r, \bn\in
   V^t\cap\zd\}$, then it is a finitely-generated $R_V$-module,
   hence Noetherian over $R_V$. Conversely, if $M$ is Noetherian
   over $R_V$ and $H\in\Hd$ contains $V$, then trivially $M$ is
   $R_H$-Noetherian as well, so that $M$ is Noetherian along~ $V$.
\end{proof}

For a subspace $V$ of $\RR^d$, let $V^\perp$ denote its
orthogonal complement. 

\begin{theorem}\label{thm:subspace-characterization}
   Let $M$ be a Noetherian $\rd$-module, $\aM$ be the
   corresponding algebraic $\zd$-action, and $V$ be a subspace of
   $\RR^d$. Then the following are equivalent.
   \begin{enumerate}
     \item[(1)] $\aM$ is expansive along~ $V$.
     \item[(2)] $\ardp$ is expansive along~ $V$ for every
      $\fp\in\asc(M)$.
     \item[(3)] $\rdp$ is Noetherian along~ $V$ and
      $V^\perp\cap\log|\V(\fp)|=\emptyset$ for every $\fp\in\asc(M)$.
   \end{enumerate}
\end{theorem}

\begin{proof}
   (1)\ifff(2): The case $V=\RR^d$ is exactly Theorem
   \ref{thm:expansive}. We may therefore assume that $\dim V\le
   d-1$. By \cite[Theorem 3.6]{BL}, $\aM$ is not expansive along~
   $V$ if and only if there is a $W\in\G_{d-1}$ containing $V$
   such that $\aM$ is not expansive along~ $W$. By Lemma
   \ref{lem:boundary-half-space}, this occurs if and only if
   there is an $H\in\N(\aM)$ containing ~$V$. Theorem
   \ref{thm:characterization} shows that this happens if and only
   if $H\in\N(\ardp)$ for some $\fp\in\asc(M)$. Reversing the
   chain of equivalences, this time for $\ardp$, shows that this
   occurs if and only if $\ardp$ is not expansive along~ $V$ for
   some $\fp\in\asc(M)$.

   (2)\ifff(3): As in the proof of (1)\ifff(2), $\ardp$ is not
   expansive along~ $V$ if and only if there is an $H\in\N(\ardp)$
   containing $V$. Theorem~ \ref{thm:characterization} shows that
   this occurs if and only if either $\rdp$ is not
   $R_H$-Noetherian or
   $[0,\infty)\vH\cap\log|\V(\fp)|\ne\emptyset$, where $\vH$ is
   the outward unit normal for~ $H$. By definition, $\ardp$ is not
   Noetherian along~ $V$ if and only if there is a half-space~ $H$
   containing~ $V$ for which $\rdp$ is not $R_H$-Noetherian. If
   $V\subset H$, then $\vH\in V^\perp$. Hence there is an
   $H\in\Hd$ containing~ $V$ with
   $[0,\infty)\vH\cap\log|\V(\fp)|\ne\emptyset$ if and only if
   $V^\perp\cap\log|\V(p)|\ne\emptyset$. 
\end{proof}

If $M$ is a Noetherian $\rd$-module and $V$ is a $k$-plane, then
Remark \ref{rem:non-Noetherian-M} and Definition
\ref{def:Noetherian-along-V} show that $M$ is Noetherian along
~$V$ if and only if $\rdp$ is Noetherian along~ $V$ for every
$\fp\in\asc(M)$. Hence
\begin{equation}\label{eqn:Nk-union}
   \Nn_k(\al_M)=\bigcup_{\fp\in\asc(M)} \Nn_k(\al_{\rdp}).
\end{equation}
Let us define
\begin{displaymath}
   \Nv_k(\ardp)=\{\,V\in\G_k:V^\perp\cap\log|\V(\fp)|\ne\emptyset\,\}
\end{displaymath}
and
\begin{displaymath}
   \Nv_k(\al_M)=\bigcup_{\fp\in\asc(M)} \Nv_k(\ardp).
\end{displaymath}
Then Theorem \ref{thm:subspace-characterization} says that
\begin{displaymath}
   \N_k^{}(\aM)=\Nn_k(\al_M)\cup\Nv_k(\al_M).
\end{displaymath}
When $k=d-1$, this is the image of the equality
\eqref{eqn:N-as-union} under the map $\pi\colon\Hd\to\G_{d-1}$
defined by $\pi(H)=\partial H$.

Next, we prove a lower-dimensional version of Proposition
\ref{prop:Nn-empty-interior}. 

\begin{proposition}\label{prop:Nn-lower-empty}
   Let $M$ be a Noetherian $\rd$-module, and let $k\le d$ be such
   that $\kdim(\rdp)\le k$ for every $\fp\in\asc(\al_M)$. Then
   $\Nn_k(\al_M)$ is a closed subset of $\G_k$ with empty interior.
\end{proposition}

\begin{proof}
   $\Nn_k(\al_M)$ is closed by Remark \ref{rem:Nnk-closed}. The case
   $d=k$ is trivial. Suppose that $k=d-1$. Let $\pi\colon
   \Hd\to\G_{d-1}$ be $\pi(H)=\partial H$. By definition,
   $\pi\bigl(\Nn(\aM)\bigr)=\Nn_{d-1}(\aM)$. Since
   $\Nn(\aM)=\bigcup_{\fp\in\asc(M)} \Nn(\ardp)$ and each
   $\Nn(\ardp)$ is nowhere dense in $\Hd$ by Proposition
   \ref{prop:Nn-empty-interior}, it follows that $\Nn_{d-1}(\aM)$ is
   nowhere dense.

   We prove the general case by downward induction on $k$. Assume
   the result is established for $k+1$. Let $\SV$ be an open set
   in $\G_k$. Then there is an open set $\SW$ in $\G_{k+1}$ such
   that every $W\in\SW$ contains a $V\in\SV$. By the induction
   hypothesis, $\SW$ contains a rational $(k+1)$-plane~ $W$ along
   which $M$ is Noetherian. By Lemma
   \ref{lem:Noetherian-along-V}, $M$ is Noetherian over
   $R_W=\ZZ[u^{\bn}:\bn\in W\cap\zd]$. We are now in the
   codimension one situation of the first part of the proof. The
   associated primes of $M$ as an $R_W$-module are $\fp\cap R_W$
   for $\fp\in\asc(M)$. Then $R_W/(\fp\cap R_W)$ is a subring of
   $\rdp$, hence the hypothesis on Krull dimension is satisfied
   for this situation. Thus there is a $V\in\SV$ along which $M$
   is Noetherian. This proves that $\Nn_k(\aM)$ is nowhere dense in
   ~$\G_k$. 
\end{proof}

\begin{example}\label{exam:better-Nn}
   We can now give another way to determine the non-Noeth\-erian set
   $\Nn$ in Example \ref{exam:eins-ward} that does not use the
   2-adic arguments there.  The polynomial $w-2\in\fp$
   shows that $\Nn$ is contained in the lower hemisphere,
   and then $1+u+v\in\fp$ shows that $\Nn$ is contained in the
   union $Q$ of the three quarter-meridians shown in Figure
   \ref{fig:eins-ward}(b).

   Assume that for some $\bv$ in $Q$ the
   module $R_3/\fp$ is Noetherian over $R_{H_{\bv}}$.
   Since the Noetherian set is open, we may assume that
   $\bv$ is rational and does not lie on the equator. Then there
   is a great circle $C$ through $\bv$ which does not contain any
   other point of $Q$. We may assume that there is an $\bn\in\zd$
   such that $C$ is the intersection of $\sph_2$ and the plane
   orthogonal to $\bn$.

   By Definition \ref{def:Noetherian-along-V}, the module
   $R_3/\fp$ is Noetherian along~ $V=\RR\bn$, since every half-space
   $\Hw$ containing $V$ has $\bw\in C$. By Lemma
   \ref{lem:Noetherian-along-V}, $R_3/\fp$ is Noetherian over the
   ring $S=\ZZ[u^{\pm\bn}]$.
   
   If $S\cap\fp\ne\{0\}$, then $\kdim(S/(\fp\cap S))=1$, and
   since $R_3/\fp$ is an integral extension we would have
   $\kdim(R_3/\fp)=1$, a contradiction. If $S\cap\fp=\{0\}$, then
   $w^{-1}=1/2\in R_3/\fp$ cannot be integral over $S$,
   contradicting the fact that $R_3/\fp$ is Noetherian over
   $S$. Hence $\Nn$ consists of all of $Q$.
\end{example}

\section{Homoclinic points and groups}\label{sec:homoclinic-group}

For the first part of this section we return to topological
$\zd$-actions. Recall from Section~ 2 that $V^t$ denotes the
thickening of a subspace $V$ by ~$t$.

\begin{definition}\label{def:homoclinic}
   Let $\beta$ be a $\zd$-action on a compact metric space
   $(X,\rho)$. Suppose that $V$ is a subspace of $\RR^d$ and that
   $y_0\in X$. We say that $x\in X$ is \textit{homoclinic to
   $y_0$ along~ $V$} if there is a $t>0$ such that
   $\rho\bigl(\al^{\bn}(x),\al^{\bn}(y_0)\bigr)\to0$ as
   $\|\bn\|\to\infty$ and $\bn\in V^t$. The set of points in $X$
   homoclinic to ~$y_0$ along ~$V$ is denoted by $\Db(y_0,V)$. In
   case $V=\RR^d$ we delete the phrase ``along~ $V$'', and write
   $\Db(y_0)$ for $\Db(y_0,\RR^d)$.
\end{definition}

\begin{remark}\label{rem:expansive-radius}
   If there is some $t_0>0$ for which
   $\rho\bigl(\al^{\bn}(x),\al^{\bn}(y_0)\bigr)\to0$ as
   $\|\bn\|\to\infty$ and $\bn\in V^{t_0}$, then this also holds
   along~ $V^t$ for every $t>0$. This follows easily from the
   observation that $\zd\cap V^{t_0}$ has bounded gaps together
   with continuity of ~$\beta$.
\end{remark}

\begin{proposition}\label{prop:countable}
   Let $V$ be an expansive subspace for the topological
   $\zd$-action ~$\beta$ on ~$X$, and $y_0\in X$. Then $\Db(y_0,V)$
   is at most a countable set.
\end{proposition}

\begin{proof}
   This can be proven exactly as in \cite[Lemma~ 3.2]{LS}.
\end{proof}

\begin{example}\label{exam:full-shift-homoclinic}
   Let $\SA$ be a finite alphabet, $X=\SA^{\zd}$, and $\beta$ be
   the $\zd$-shift action on $X$. Then $\beta$ is expansive, and
   $\Db(y_0)$ is the countable set of points in $X$ differing
   from $y_0$ in only finitely many coordinates. For every
   subspace $V$ of dimension $<d$, it is easy to see that
   $\Db(y_0,V)$ is uncountable.
\end{example}

\begin{example}\label{exam:ledrappier-homoclinic}
   We return to Ledrappier's example (see Examples \ref{exam:ledrappier}
   and \ref{exam:ledrappier-set}). Here
   $\beta=\al_{R_2/\< 2,1+u+v\>}$ and $X=X_{R_2/\<
     2,1+u+v\>}$. We choose $y_0=0_X=0$. Since the only
   point in $X$ having finitely many nonzero coordinates is $0$,
   we see  that
   $\Db(0)=\{0\}$. Let $L_{\theta}$ denote the line in $\RR^2$
   making angle $\theta$ with the positive horizontal axis. For
   all $\theta$ with $0<\theta<\pi/2$, the sets
   $\Db(0,L_{\theta})$ are all equal. Each consists of points of
   the form indicated in Figure
   \ref{fig:ledrappier-homoclinic}(a), where the shaded regions
   contain only 0's, and coordinates represented by the dots can
   be of an arbitrary finite length, filled in with arbitrary
   values, and these are used to determine all remaining
   coordinates.

   \begin{figure}[htbp]
      \begin{center}
          \scalebox{1.0}{\includegraphics{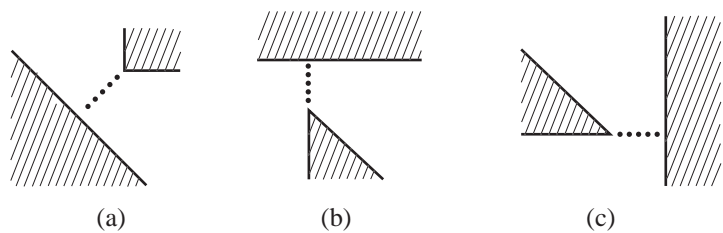}}
         \caption{Homoclinic points for Ledrappier's example}
         \label{fig:ledrappier-homoclinic}
      \end{center}
   \end{figure}

   For $\pi/2<\theta<3\pi/4$, all the sets $\Db(0,L_{\theta})$
   are also equal, and each consists of points shown in Figure
   \ref{fig:ledrappier-homoclinic}(b), with the same conventions
   as before. For $3\pi/4<\theta<\pi$, Figure
   \ref{fig:ledrappier-homoclinic}(c) describes the form of
   points in each $\Db(0,L_{\theta})$. It is also easy to see
   that $\Db(0,L_0)=\Db(0,L_{\pi/2})=\Db(0,L_{3\pi/4})=\{0\}$.

   Here $\beta$ has expansive components
   $\SC_1=\{L_{\theta}:0<\theta<\pi/2\}$,
   $\SC_2=\{L_{\theta}:\pi/2<\theta<3\pi/4\}$, and
   $\SC_3=\{L_{\theta}:3\pi/4<\theta<\pi\}$. Within an expansive
   component the homoclinic set is constant, but it changes
   abruptly when passing from one component to another. See
   \cite{MS} for another example of how the homoclinic group varies.
\end{example}

The following result shows that homoclinic sets are always
constant within an expansive component. Roughly speaking, nearby
planes in an expansive component code each other, so a point that
is homoclinic for one is also homoclinic for the other.

\begin{theorem}\label{thm:homoclinic-constancy}
   Let $\beta$ be a topological $\zd$-action on a compact metric
   space~ $X$, let $\SC$ be a connected component of\,
   $\E_k(\beta)$, and let $y_0\in X$. Then
   $\Db(y_0,V)=\Db(y_0,W)$ for all $V,W\in\SC$.
\end{theorem}

\begin{proof}
   We show that $\Db(y_0,\cdot)$ is locally constant on ~$\SC$,
   hence constant on ~$\SC$ by connectedness.
   
   Let $U\in\SC$ and $\delta$ be an expansive constant for
   ~$\beta$. Recall Definition \ref{def:coding} of coding for
   subsets of ~$\RR^d$. Let $B(r)$ be the open ball of radius $r$
   in $\RR^d$. For $V\in\G_k$ and $\eta>0$ put
   $V_{\eta}=\{\bx\in\RR^d :\dist(\bx,V)<\eta\|\bx\|\}$. It follows
   from \cite[Proposition~ 3.8]{BL} that there is a compact
   neighborhood $\SU$ of ~$U$ in ~$\SC$ and positive numbers
   $s,t,\eta>0$ such that for every $V,W\in\SU$ and every $r>0$
   we have that $V^t\smallsetminus B(r)$ codes
   $W_\eta\smallsetminus B(r+s)$.

   Suppose that $V,W\in\SU$ and $x\in\Db(y_0,V)$. Then there is
   an $r>0$ such that
   $\rho\bigl(\beta^{\bn}(x),\beta^{\bn}(y_0)\bigr)\le\delta$ for
   every $\bn\in V^t\smallsetminus B(r)$. Hence by coding, this
   inequality also holds for every $\bn\in W_\eta\smallsetminus
   B(r+s)$. Let  $\varepsilon>0$. A simple compactness argument
   using expansiveness of $\beta$ shows that there is an $a>0$
   such that $\rho_\beta^{B(a)}(x,y_0)\le\delta$ implies that
   $\rho(x,y_0)<\varepsilon$. Choose $q$ large enough so that if
   $\bx\in W$ and $\|\bx\|>q$, then $\bx+B(a+t)\subset
   W_\eta\smallsetminus B(r+s)$. It follows that
   $\rho\bigl(\beta^{\bn}(x),\beta^{\bn}(y_0)\bigr)<\varepsilon$
   for every $\bn \in W^t$ with $\|\bn\|>q$. This proves that
   $x\in\Db(y_0,W)$, and so $\Db(y_0,V)\subset\Db(y_0,W)$.
   Interchanging the roles of $V$ and $W$ gives the reverse
   inclusion. 
\end{proof}

We next turn to considering homoclinic points for an algebraic
$\zd$-action $\al$ on ~$X$. Since $\rho$ is
translation invariant,
$\rho\bigl(\al^{\bn}(x),\al^{\bn}(y_0)\bigr)=
\rho\bigl(\al^{\bn}(x-y_0),0_X\bigr)$, so that
$\Da(y_0,V)=y_0+\Da(0_X,V)$. Hence for algebraic $\zd$-actions we
need only compute $\Da(0_X,V)$, which we shorten to $\Da(V)$.
Obviously $\Da(V)$ is a subgroup of ~$X$, which we term the
\textit{homoclinic group of $\al$ along ~$V$}. When $V=\RR^d$, the
subgroup $\Da=\Da(\RR^d)$ is called the \textit{homoclinic group
of ~$\al$}.

The homoclinic groups of algebraic $\zd$-actions were studied in
\cite{LS}, especially for expansive actions. The main result
there is the following \cite[Theorems 4.1 and 4.2]{LS}. Here
``entropy'' is taken with respect to Haar measure.

\begin{theorem}\label{thm:homoclinic-group}
   Let $\al$ be an expansive algebraic $\zd$-action on ~$X$.
   \begin{enumerate}
     \item[(1)] $\Da$ is nontrivial if and only if $\al$ has positive
      entropy.
     \item[(2)] $\Da$ is nontrivial and dense in $X$ if and only
      if $\al$ has completely positive entropy.
   \end{enumerate}
\end{theorem}

Using this together with constancy of the homoclinic group within
an expansive component, we obtain further instances of the
expansive subdynamics philosophy, as follows.

Call $V\in\G_k$ \textit{rational} if $V\cap\zd$ spans $V$. For an
expansive component $\SC\subset\E_k(\al)$ let
$\SCQ=\{V\in\SC:\text{ V is rational}\}$, which is a dense
subset of $\SC$. For $V\in\SCQ$ we let $\al|_{V\cap\zd}$ denote
the rank $k$ action obtained by restricting $\al$ to $V\cap\zd$. 

\begin{theorem}\label{thm:entropy}
   Let $\al$ be an algebraic $\zd$-action and $\SC$ be an
   expansive component of $\E_k(\al)$. If for some $V\in\SCQ$ the
   action $\al|_{V\cap\zd}$ has positive entropy, has completely
   positive entropy, or is isomorphic to a $\ZZ^k$ Bernoulli
   shift, then the same property holds for $\al|_{W\cap\zd}$
   for every $W\in\SCQ$.
\end{theorem}

\begin{proof}
   By Theorem~ \ref{thm:homoclinic-group}, $\al|_{V\cap\zd}$ has
   positive entropy if and only if $\Da(V)\ne\{0_X\}$, and by Theorem~
   \ref{thm:homoclinic-constancy} we have $\Da(V)=\Da(W)$ for all
   $W\in\SCQ$, establishing the positive entropy portion. The
   proof for completely positive entropy is similar, using
   density rather than nontriviality of the homoclinic group.
   Finally, completely positive entropy for an algebraic
   $\zd$-action is equivalent to Bernoullicity (see
   \cite{RS} or \cite[\S23]{Sch}).
\end{proof}

\begin{remarks}
   (1) By \cite[Theorem~ 6.3(4)]{BL}, if there is a rational
   $V\in\E_k(\al)$ for which $\al|_{V\cap\zd}$ has zero entropy,
   then $\al|_{W\cap\zd}$ has zero entropy for every $W\in\G_k$.
   This not only provides an alternative proof of the first part
   of the previous theorem, it also shows that on the entire set
   $\E_k(\al)$ either entropy vanishes identically or it is
   strictly positive everywhere.

   (2) It is likely that suitable definitions of completely
   positive entropy and Bernoullicity can be found for
   non-rational $k$-planes, and that the conclusions of Theorem~
   \ref{thm:entropy} can be extended to all $W\in\SC$. However,
   we will not pursue this further here.
\end{remarks}

\begin{example}\label{exam:hc-groups-for-planes}
   (\textit{Homoclinic groups of $2$-planes for Example
   \ref{exam:eins-ward}}.) To illustrate some of the phenomena
   surrounding the preceding theorems, we describe the homoclinic
   groups of $2$-planes for Example \ref{exam:eins-ward}. In this
   example, $d=3$, $\fp=\<1+u+v,w-2\>$, $X=X_{\rthp}$, and
   $\al=\arthp$. Observe that since $\fp$ is nonprincipal, $\al$
   has zero entropy. Since $\al$ is expansive, $\Da=\{0\}$ by Theorem
   \ref{thm:homoclinic-group}.
   
   To help describe the groups $\Da(V)$ for $2$-planes
   $V\in\G_2$, recall the map $\pi\colon \H_3\to\G_2$ defined by
   $\pi(H)=\partial H$. Under our standing correspondence
   $\H_3\leftrightarrow \sph_2$, this map is the usual
   identification of antipodal points of $\sph_2$ to obtain
   projective $2$-space $\G_2$. By vertically projecting the
   upper hemisphere of $\sph_2$ to the unit disk $\D$, we can
   represent $\G_2$ as $\D$ with antipodal boundary points
   identified. Using the representation, the image $\pi(\N(\al))$
   in $\D$ is shown in Figure \ref{fig:two-planes}(a). The shaded
   region is $\pi(\Nv(\al))$, while the three segments comprise
   $\pi(\Nn(\al))$. There are three expansive components of
   $2$-planes, labeled $\SC_1$, $\SC_2$, and $\SC_3$.

   \begin{figure}[htbp]
      \begin{center}
         \scalebox{.85}{\includegraphics{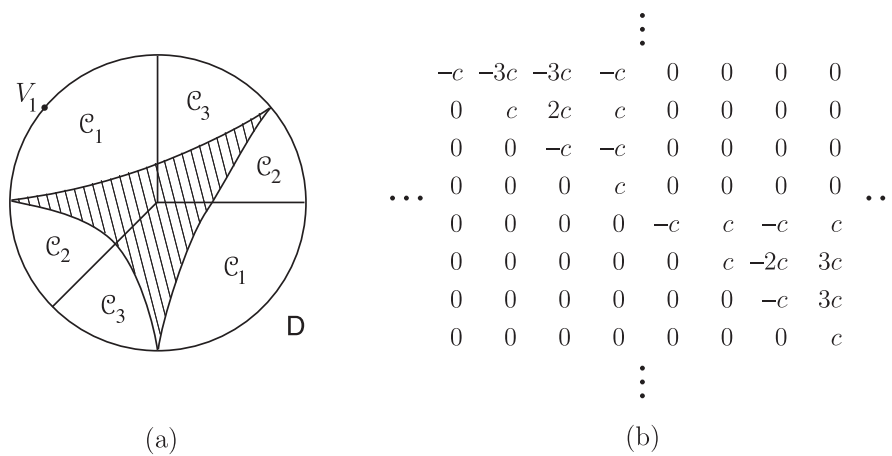}}
         \caption{Analysis for homoclinic groups of $2$-planes}
         \label{fig:two-planes}
      \end{center}
   \end{figure}

   We say that $y\in\Da(V)$ is a \textit{fundamental homoclinic
   point for $V$} if $\{\al^{\bn}(y):\bn\in\ZZ^3\}$ generates
   $\Da(V)$ as an abelian group.

   We will sketch the following description of how $\Da(V)$
   varies for $V\in\D$.
   \begin{enumerate}
     \item[(1)] For each expansive component $\SC_r$ there is an
      explicit nontrivial fundamental homoclinic point $x^r$
      (recall that $\Da(\,\cdot\,)$ is constant within an expansive
      component),
     \item[(2)] $\Da(V)=\{0\}$ for every $V\in\pi(\Nn(\al))$, and
     \item[(3)] For every
      $V\in\pi(\Nv(\al))\smallsetminus\pi(\Nn(\al))$ there is a
      nontrivial fundamental homoclinic point $x^{\!V}$ that varies
      continuously with $V$.
   \end{enumerate}
   
   For (1), first consider component $\SC_1$. Define $x^1\in X$
   by
   \begin{displaymath}
      x_{i,j,k}^1=
      \begin{cases}
         (-1)^{i-1}\dbinom{i}{-j}\,2^k   &\text{if $i\ge0$,
         $j\le0$,} \\[10pt]
         (-1)^{j-1}\dbinom{j-1}{-i-1}\,2^k &\text{if $i\le -1$,
         $j\ge1$,} \\
         0         &\text{otherwise.}
      \end{cases}
   \end{displaymath}
   Figure \ref{fig:two-planes}(b) shows a typical horizontal
   slice of $x^1$ at level $k$, where $c=2^k$. The $2$-plane
   $V_1\in\SC_1$, shown in Figure \ref{fig:two-planes}(a), is the
   vertical plane that intersects the horizontal $uv$-plane in
   the line $u=v$.  Clearly $x^1\in\Da(V_1)$. Using the fact that
   the coordinates of any point in $X$ homoclinic along a
   vertical line must be an integer times the successive powers
   of $2$, it is not hard to show that integer combinations of
   shifts of $x^1$ comprise all of $\Da(V_1)$. Thus $x^1$ is a
   fundamental homoclinic point for $V_1$. It is therefore also a
   fundamental homoclinic point for all $V\in\SC_1$. A similar
   construction produces fundamental homoclinic points $x^2$ and
   $x^3$ for $\SC_2$ and $\SC_3$, respectively, which bear the
   same relationship to $x^1$ as Figures
   \ref{fig:ledrappier-homoclinic}(b) and
   \ref{fig:ledrappier-homoclinic}(c) bear to Figure
   \ref{fig:ledrappier-homoclinic}(a).
   
   It is also possible to show using this description that the
   intersection of the homoclinic groups of any pair of distinct
   expansive components is trivial, analogous to a result in
   \cite{MS} for a pair of commuting toral automorphisms.

   For (2), we give a complete argument only for the horizontal
   plane $V_0$. A somewhat more complicated proof is required for
   general $V\in\pi(\Nn(\al))$. For brevity we assume familiarity
   with \cite{LS}.

   Suppose that $x\in\Da(V_0)$. Let $f(u,v)=1+u+v$, and define
   $\ftil\in\ell^\infty(\ZZ^3,\ZZ)$ by
   $\ftil_{\bn}=c_f(-\bn)$. Choose
   $\xbar\in\ell^\infty(\ZZ^3,\RR)$ so that
   $|\xbar_{\bn}|\le\frac12$ and $\xbar_{\bn}\equiv
   x_{\bn}\pmod1$. Then $\xbar_{\bn}\to0$ along each horizontal
   plane $V_k=V_0+k\be_3$. Since $x\in X$, we have that
   $\ftil*\xbar\in\ell^\infty(\ZZ^3,\ZZ)$ and also tends to $0$
   in each $V_k$. Hence there are polynomials $h_k\in
   R_2/\<1+u+v\>$ such that
   $\bigl. \ftil*\xbar\bigr|_{V_k}=(w^kh_k)\sptilde$. Since
   $w-2\in\fp$, it follows that $2h_{k-1}\equiv h_k$ in
   $R_2/\<1+u+v\>$. Hence each $h_k$ is divisible in
   $R_2/\<1+u+v\>$ by arbitrarily large powers of $2$. But
   $R_2/\<1+u+v\>\cong \ZZ[u^{\pm1},1/(u+1)]$ is a localization
   of $\ZZ[u]$ in which $2$ is not invertible. Hence each
   $h_k=0$. Let $\ybar\in\ell^\infty(\ZZ^2,\RR)$ be the
   restriction of $\xbar$ to some $V_k$. Then $\ftil*\ybar=0$, so
   that $\ybar$ is a pseudo-measure on $\ZZ^2$ whose support must
   be contained in the finite set $\V(f)\cap\SS^2$. But then
   $\ybar$ is almost periodic, and in particular $\ybar\notin
   c_0(\ZZ^2)$ unless $\ybar=0$. This proves that the only point
   homoclinic for $V_0$ is $x=0$.
   
   For (3), let $V\in\pi(\Nv(\al))\smallsetminus\pi(\Nn(\al))$.
   Choose $a,b$ so that $(a,b,1)$ is normal to $V$. Put
   $f_{a,b}(u,v)=1+2^au^{-1}+2^bv^{-1}$ and $F_{a,b}=1/f_{a,b}$. To
   construct $x^{\!V}$, we first construct an auxiliary point
   $y\in\ell^\infty(\ZZ^2,\RR)$ as follows. If $V\notin\partial
   \pi(\Nv(\al))$ then $F_{a,b}$ has two poles of order one in
   $\SS^2$, and so $F_{a,b}\in L^1(\SS^2)$. In this case we put
   $y_{\bn}=\widehat{F}_{a,b}(\bn)$. Then
   \begin{equation}\label{eqn:riemann-lebesgue}
      y_{\bn}\to0 \quad\text{as\quad $\|\bn\|\to\infty$}
   \end{equation}
   by the Riemann-Lebesgue Lemma. Also,
   \begin{equation}\label{eqn:2^a-ledrappier}
      (\ftil_{a,b}*y)_{\bn} =
      y_{\bn}+2^ay_{\bn+\be_1}+2^by_{\bn+\be_2} =
      \begin{cases}
         1 &\text{if $\bn=\mathbf{0}$},\\
         0 &\text{otherwise}.
      \end{cases}
   \end{equation}
   If $V\in\partial\pi(\Nv(\al))$, we can still expand $F_{a,b}$
   in a series. For example, if $a<0$ and $b<0$, then $V$ being
   on the boundary of $\pi(\Nv(\al))$ corresponds to
   $2^a+2^b=1$. Hence
   \begin{align*}
      F_{a,b}&=\frac1{1+2^au^{-1}+2^bv^{-1}} =
      \sum_{k=0}^\infty (-1)^k(2^au^{-1}+2^bv^{-1})^k \\
      &=\sum_{k=0}^\infty\sum_{j=0}^k (-1)^k\dbinom{k}{j}
      2^{aj+b(k-j)}u^{-j}v^{-(k-j)}.
   \end{align*}
   In this case let $y_{(m,n)}$ denote the coefficient of
   $u^mv^n$ in the expansion. Then \eqref{eqn:riemann-lebesgue}
   holds by standard estimates and \eqref{eqn:2^a-ledrappier} by
   construction.

   Next we construct $x^{\!V}$ from $y$. Let
   $V'=V-[0,1)\be_3$. Vertical projection from $\RR^3$ to $\RR^2$
   gives a bijection from $V'\cap\ZZ^3$ to $\ZZ^2$. Let $\{r\}$
   denote the fractional part of a real number $r$. Define
   $x_{\bn}^{\!V}$ for $\bn \in V'\cap\ZZ^3$ by
   \begin{displaymath}
      x^{\!V}_{(m,n,\lfloor -ma-nb\rfloor)}=2^{-\{-ma-nb\}}y_{(m,n)}.
   \end{displaymath}
   Extend $x^{\!V}$ to those $\bn$ above $V'$ using
   $x^{\!V}_{\bn+\be_3}=2x^{\!V}_{\bn}$, and to those $\bn$ below $V'$
   using $x^{\!V}_{\bn}+x^{\!V}_{\bn+\be_1}+x^{\!V}_{\bn+\be_2}=0$. The
   former extension is clearly unique, and the latter is unique
   since $V\notin\pi(\Nn(\al))$. An elementary argument shows
   that $x^{\!V}\in X$, and $x^{\!V}$ is homoclinic along $V$ by
   construction. It is not difficult to establish that $x^{\!V}$ is a
   fundamental homoclinic point for $V$. Finally, the
   construction shows that $x^{\!V}$ varies continuously for
   $V\in\pi(\Nv(\al))\smallsetminus \pi(\Nn(\al))$.
\end{example}

\end{document}
